\pdfoutput=1
\documentclass[12pt]{amsart}
\usepackage[T1]{fontenc}
\usepackage{amsmath,amssymb,mathtools}
\usepackage{geometry}
\usepackage{float}
\usepackage{needspace}
\usepackage{pgfplots}
\pgfplotsset{compat=1.18}
\usepackage[unicode,hidelinks]{hyperref}
\usepackage[alphabetic,msc-links]{amsrefs}
\hypersetup{pdftitle={Noncompactifiable Kahler threefolds with positive Ricci curvature},pdfauthor={Gang Liu and Bo Zhu}}
\allowdisplaybreaks[2]
\numberwithin{equation}{section}
\newtheorem{theorem}{Theorem}[section]
\newtheorem{proposition}[theorem]{Proposition}
\newtheorem{lemma}[theorem]{Lemma}
\newtheorem{corollary}[theorem]{Corollary}
\newtheorem{conjecture}[theorem]{Conjecture}
\theoremstyle{definition}
\newtheorem{definition}[theorem]{Definition}
\newtheorem{remark}[theorem]{Remark}
\newcommand{\C}{\mathbb C}
\newcommand{\R}{\mathbb R}
\newcommand{\Z}{\mathbb Z}
\newcommand{\QQ}{\mathbb Q}
\newcommand{\ii}{\sqrt{-1}}
\newcommand{\dd}{\mathrm d}
\newcommand{\ddc}{\ii\partial\bar\partial}
\newcommand{\eps}{\varepsilon}
\newcommand{\Om}{\Omega}
\newcommand{\wtO}[1]{O_{\mathrm{wt}}\!\left(#1\right)}
\DeclareMathOperator{\Ric}{Ric}
\DeclareMathOperator{\Vol}{Vol}

\DeclareMathOperator{\Res}{Res}
\DeclareMathOperator{\SO}{SO}
\DeclareMathOperator{\arsinh}{arsinh}
\DeclareMathOperator{\arcosh}{arcosh}
\DeclareMathOperator{\id}{id}

\title[Noncompactifiable K\"ahler threefolds]{Noncompactifiable K\"ahler threefolds with positive Ricci curvature}
\author{Gang Liu}
\address[Gang Liu]{School of Mathematical Sciences, Shanghai Key Laboratory of PMMP, East China Normal University}
\email{gliu@math.ecnu.edu.cn}
\thanks{Gang Liu is partially supported by NSFC No.~12071140.}
\author{Bo Zhu}
\address[Bo Zhu]{Yau Mathematical Sciences Center, Tsinghua University}
\email{zhub@tsinghua.edu.cn}
\thanks{Bo Zhu is supported by NSFC No.~12501066.}
\date{November 14, 2026}

\begin{document}
\begin{abstract}
We construct a complete K\"ahler threefold \((M,\omega)\) with positive Ricci curvature
and infinite-dimensional third rational homology. It admits no analytic
compactification, answering Yau's compactification question negatively
(see \cite{Yau1992}{}, Problem~71, p.~304). The metric is obtained by solving a complex
Monge--Amp\`ere equation on an affine threefold and lifting the solution
to an infinite cyclic cover. Its Ricci form \(\rho_\omega\) satisfies
\(\int_M\rho_\omega^3=\infty\), disproving the top Yang--Yau conjecture.
\end{abstract}
\maketitle

\section{Introduction}\label{sec:introduction}

The Frankel conjecture, proved by Mori and by Siu--Yau, states that
a compact K\"ahler
\(n\)-manifold with positive holomorphic bisectional curvature is
biholomorphic to \(\mathbb P^n\) (see \cites{Mori,SiuYau}).

A compact K\"ahler manifold admits a metric of positive Ricci curvature
exactly when it is Fano, with the converse following from Yau's solution
of the Calabi conjecture~\cite{YauCalabi}. Such a manifold is simply
connected by Kobayashi's theorem~\cite{KobayashiRicci}.

For complete noncompact K\"ahler manifolds, S.-T.~Yau proposed the
following compactification conjecture in \emph{Open problems in geometry}
(see \cite{Yau1992}{}, Problem~71, p.~304):
\begin{conjecture}[Yau's compactification conjecture]\label{conj:yau-compactification}
Every complete K\"ahler manifold \((M,\omega)\) with positive Ricci
curvature is biholomorphic to a Zariski open subset of a compact
K\"ahler manifold.
\end{conjecture}
Here, ``Zariski open'' means the complement of a closed complex analytic
subset of the ambient manifold.
A stronger assertion is that \(M\) is biholomorphic to a
quasi-projective variety.

N.~Mok proved that a complete K\"ahler manifold \((M^n,\omega)\)
with positive Ricci curvature is quasi-projective if it has maximal
volume growth, quadratic curvature decay, and
\(\int_M\rho_\omega^n<\infty\), where \(\rho_\omega\) is the Ricci form
(see \cite{Mok}{}, Main Theorem, p.~373).
W.-K.~To obtained quasi-projective embeddings without assuming maximal
volume growth. His quadratic scalar-curvature decay case requires
positive Ricci curvature, finite topological type, finite
top Ricci integral, and
\(\lvert\operatorname{Bisect}(v,w)\rvert\le C\Ric(v,\bar v)\)
for unit \((1,0)\) vectors \(v,w\) (see \cite{To}{}, main theorem).
The first author proved quasi-projectivity without the finiteness
assumption on the top Ricci integral, under positive Ricci curvature,
maximal volume growth, and quadratic decay of the full curvature tensor
(see \cite{Liu}{}, Theorem~1.2(III)).

To generalize the Cohn--Vossen inequality to higher dimensions, Yau
proposed the following conjecture (see \cite{Yau1992}{}, Problem~9).
\begin{conjecture}[Yau's higher-dimensional Cohn--Vossen conjecture]\label{conj:yau-cohn-vossen}
Let \((M,g)\) be a complete noncompact Riemannian manifold of real
dimension \(m\) with nonnegative Ricci curvature. Then, for every
\(p\in M\),
\[
 \limsup_{R\to\infty}R^{2k-m}
 \int_{B(p,R)}\sigma_k(\Ric)\,\dd V<\infty,
 \qquad 1\le k\le m,
\]
where \(\sigma_k(\Ric)\) is the \(k\)-th elementary symmetric
function of the eigenvalues of the real Ricci tensor.
\end{conjecture}
Bo Yang constructed complete \(U(n)\)-invariant K\"ahler metrics on
\(\C^n\) with nonnegative holomorphic bisectional curvature that
disprove these estimates for \(2\le k\le n\), for every \(n\ge2\)
\cite{YangCohnVossen}.
Hao and J.~Zhu announced a counterexample to the
\(k=1\) scalar-curvature case on a complete Riemannian
three-manifold with strictly positive Ricci curvature
\cite{HaoZhu}.
Xu also constructed a sequence of smooth metrics \(g_j\) on \(S^3\)
with \(\Ric_{g_j}>0\), \(\operatorname{diam}(S^3,g_j)<1/2\), and
\[
 \int_{S^3}\operatorname{Scal}_{g_j}\,\dd V_{g_j}\longrightarrow\infty
 \qquad\text{as }j\to\infty
\]
(see \cite{XuCollapse}, Theorem~1.1).
Since every unit ball is the whole manifold, these compact examples
disprove the corresponding dimension-only local scalar-curvature
integral bound under nonnegative Ricci curvature.
The sequence is volume-collapsing.

In the K\"ahler setting, Yang proposed an analogue using wedge powers
of the Ricci form (see \cite{YangCohnVossen}, Question~1.2).
We refer to the following formulation under nonnegative Ricci
curvature as the Yang--Yau conjecture
(cf. \cite{DatarPingaliSeshadri}).
\begin{conjecture}[Yang--Yau]\label{conj:yang-yau}
Let \((M,\omega)\) be a complete noncompact K\"ahler \(n\)-manifold
whose Ricci form satisfies \(\rho_\omega\ge0\). Then, for every
\(p\in M\),
\begin{equation}\label{eq:kahler-cohn-vossen}
 \limsup_{R\to\infty}R^{2k-2n}
 \int_{B_\omega(p,R)}\rho_\omega^k\wedge\omega^{n-k}<\infty,
 \qquad 1\le k\le n.
\end{equation}
\end{conjecture}
The top-degree case \(k=n\), called the top Yang--Yau conjecture,
asserts that
\(\int_M\rho_\omega^n<\infty\).
The first author proved this finiteness under nonnegative holomorphic
bisectional curvature and Euclidean volume growth
\cite{LiuCohnVossen}.
V.~Datar, V.~P.~Pingali, and H.~Seshadri proved finiteness
under positive holomorphic bisectional curvature
(see \cite{DatarPingaliSeshadri}{}, Theorem~1.1).

B.-L.~Chen and X.-P.~Zhu used finiteness of the top Ricci integral to
prove quasi-projectivity for complete noncompact K\"ahler manifolds
with bounded positive sectional
curvature (see \cite{ChenZhu}{}, Main Theorem).
Datar, Pingali, and Seshadri combined their finiteness result with
Chen--Zhu's theorem to obtain quasi-projectivity under bounded positive
sectional curvature (see \cite{DatarPingaliSeshadriSectional}{}, Theorem~1.2).

We show that positive Ricci curvature alone does not imply
compactifiability, even in complex dimension three. The example is
an infinite cyclic cover of a smooth affine threefold.
Its top Ricci integral also diverges: \(\int_M\rho_{\omega_M}^{\,3}=\infty\).
Thus both Conjecture~\ref{conj:yau-compactification} and the top-degree
case of Conjecture~\ref{conj:yang-yau} fail in complex dimension three.

We consider the affine threefold
\begin{equation}\label{eq:space}
 X=\left\{(u,z)\in\C^*\times\C^3:
 Q(z):=\sum_{j=1}^3z_j^2=P(u):=u+u^{-1}\right\}.
\end{equation}
Let
\[
 \Om=\Res_{Q-P=0}
 \frac{\dd\log u\wedge\dd z_1\wedge\dd z_2\wedge\dd z_3}{Q-P}
\]
be its nonvanishing holomorphic residue form.

\begin{theorem}\label{thm:main}
There exists a smooth complete K\"ahler form \(\omega\) on \(X\) such that
\begin{equation}\label{eq:mainvolume}
 \omega^3=
 \frac{6(\ii)^9\Om\wedge\bar\Om}
 {\sqrt{(1+|z|^2)^2+|u-u^{-1}|^2}}.
\end{equation}
Its Ricci form satisfies
\[
 \rho_\omega=\frac12\ddc\log
 \bigl((1+|z|^2)^2+|u-u^{-1}|^2\bigr)>0.
\]
Moreover, \(b_1(X)\ge1\) and \(\pi_1(X)\) surjects onto \(\Z\).
\end{theorem}

Theorem~\ref{thm:main} and the topology of the following cover
give the counterexample.

\begin{corollary}\label{thm:cover}
The complex threefold
\[
 M=\left\{(Z,w)\in\C^3\times\C:
 \sum_{j=1}^3Z_j^2=\sin w\right\}
\]
is a connected infinite cyclic holomorphic cover of \(X\), with covering
map \(\pi(Z,w)=(-\ii e^{\ii w},\sqrt2 Z)\).
For a form \(\omega\) provided by Theorem~\ref{thm:main}, the pullback
\(\omega_M=\pi^*\omega\) is complete and has strictly
positive Ricci curvature. We also have
\[
 \dim_\QQ H_3(M;\QQ)=\infty,
 \qquad \int_M\rho_{\omega_M}^3=\infty.
\]
In particular, \(M\) is not biholomorphic to an analytic Zariski open
subset of any compact complex manifold.
\end{corollary}

Products give further counterexamples to
Conjectures~\ref{conj:yau-compactification} and~\ref{conj:yang-yau}
in every higher dimension.

\begin{corollary}\label{cor:ricci-integrals}
For every pair of integers \(3\le k\le n\), there is a complete
noncompact K\"ahler \(n\)-manifold
\((N_{n,k},\omega_{n,k})\) with positive Ricci curvature such that,
for every \(p\in N_{n,k}\),
\[
 \lim_{R\to\infty}R^{2k-2n}
 \int_{B_{\omega_{n,k}}(p,R)}
       \rho_{\omega_{n,k}}^k\wedge\omega_{n,k}^{n-k}
 =\infty.
\]
Moreover, \(N_{n,k}\) admits no analytic compactification.
In particular, the finiteness assertion
\(\int_N\rho_\omega^n<\infty\) fails in every complex dimension
\(n\ge3\).
\end{corollary}

Tian--Yau constructed complete K\"ahler metrics with prescribed
Ricci form on suitable quasi-projective manifolds
\cites{TianYauI,TianYauII}, including positively Ricci-curved metrics
(\cite{TianYauI}, Theorem~4.3). We apply Hein's existence theorem
to the noncompact Monge--Amp\`ere equation.

\medskip
\noindent\textbf{Proof strategy.}
Let \(\nu\) denote the volume form in \eqref{eq:mainvolume}.
Section~\ref{sec:reduction} states Hein's theorem in Y.~Li's formulation
and lists its hypotheses as Conditions A, B, and C.
Section~\ref{sec:background} constructs the complete background metric
\(g_0\) by gluing a homogeneous outer potential to a smooth corrected
fibre potential and extending the result over the compact core.
Appendix~\ref{app:model} constructs the potentials and estimates the
volume error. In Section~\ref{sec:analysis}, the uniform holomorphic
charts complete the verification of A and give Condition B with
decay exponent \(21/10\). Appendix~\ref{app:atlas} proves the chart
estimates. Annular connectivity, volume bounds, and a Ricci lower
bound give \(SOB(12/5)\), which is Condition C.
Hein's theorem then gives a complete K\"ahler form satisfying
\(\omega^3=\nu\), so its Ricci form is the positive form \(\rho_\nu\).

In Section~\ref{sec:topology}, the fibre two-spheres over consecutive
zeros of \(\sin w\) close up to form three-spheres on the cyclic
cover. Their intersection pairings detect infinitely many independent
classes in \(H_3(M;\QQ)\). This excludes compactification, since an
analytic Zariski open subset of a compact complex manifold has finite
CW type. Disjoint deck translates of a lifted ball force the top Ricci
integral to diverge; products then give the higher-dimensional examples.

\medskip
\noindent\textbf{Notation.}\par
\noindent
Throughout, \(C>0\) denotes a constant that may change from line to line.
We write \(A\asymp B\) if \(C^{-1}B\le A\le CB\) for a constant
\(C\ge1\) uniform on the region in question.

\subsection*{Acknowledgements}
The authors are grateful to their former advisor, Professor J.~Wang,
and to Professor H.-J.~Hein for their interest in this work and
their helpful suggestions.

\subsection*{Declaration on AI use}
The authors led the research, using AI to assist with language editing,
literature searches, and mathematical discussion. This included reviewing
the gluing estimates, the application of Hein's theorem, and the
topological arguments. The authors take full responsibility for the final
text and its mathematical content.

\section{Reduction to Hein's theorem}\label{sec:reduction}

\begin{lemma}\label{lem:geometry}
The hypersurface \(X\) is a smooth, path-connected complex manifold,
and the function \(t:X\to[0,\infty)\), \(t(u,z)=|z|^2\), is proper.
The residue form \(\Om\) is a nowhere-vanishing holomorphic
section of \(K_X\).
\end{lemma}
\begin{proof}
The defining function of \(X\) satisfies
\[
 \dd(Q-P)=-(u-u^{-1})\,\frac{\dd u}{u}
 +2\sum_{j=1}^3z_j\,\dd z_j.
\]
If this differential vanished at a point of \(X\), then \(z=0\)
and \(u-u^{-1}=0\). This would give \(u=\pm1\), \(Q(z)=0\), and
\(P(u)=\pm2\), contradicting \(P=Q\). Thus \(X\) is smooth.
The form \(\dd\log u\wedge\dd z_1\wedge\dd z_2\wedge\dd z_3\)
is a global nowhere-vanishing holomorphic volume form on
\(\C^*\times\C^3\). Since \(\dd(Q-P)\ne0\) along \(X\), the
residue form \(\Om\) is holomorphic and nowhere vanishing. For example,
on the chart \(z_3\ne0\),
\[
 \Om=\frac{\dd\log u\wedge\dd z_1\wedge\dd z_2}{2z_3}.
\]
The analogous expressions on the other charts give the same global
form.

For properness,
\(|u+u^{-1}|=|Q(z)|\le t\le T\) bounds both \(|u|\) and
\(|u|^{-1}\) by the quadratic formula. Together with \(|z|^2=t\le T\),
this shows that the closed sublevel set \(\{t\le T\}\) is compact.

The projection \(\pi_0:X\to\C^*\), \(\pi_0(u,z)=u\),
has smooth connected affine quadrics as its generic fibres.
It is locally trivial over the connected set \(\{P\ne0\}\), so the
union of these fibres is connected and dense in \(X\).
Thus \(X\) is connected, hence path connected since it is smooth.
\end{proof}

We define the prescribed volume form by
\begin{equation}\label{eq:nu}
 \nu=
 \frac{6(\ii)^9\Om\wedge\bar\Om}
 {\sqrt{(1+t)^2+|u-u^{-1}|^2}},\qquad t=|z|^2.
\end{equation}
We seek a K\"ahler form \(\omega\) with \(\omega^3=\nu\)
whose associated metric is complete.

\smallskip
In local holomorphic coordinates \((w_1,w_2,w_3)\), we write
\[
 \omega=\ii g_{j\bar k}\dd w_j\wedge\dd\bar w_k,\qquad
 \Om=h\,\dd w_1\wedge\dd w_2\wedge\dd w_3,
\]
where \(h\) is holomorphic and nowhere zero. We use the conventions
\[
 g=2\operatorname{Re}(g_{j\bar k}\dd w_j\otimes\dd\bar w_k),
 \qquad \dd V_g=\frac{\omega^3}{6}.
\]
Since
\[
 \frac{\omega^3}{(\ii)^9\Om\wedge\bar\Om}
 =\frac{6\det(g_{j\bar k})}{|h|^2},\qquad
 \ddc\log|h|^2=0,
\]
the Ricci form is
\[
 \rho_\omega=-\ddc\log\det(g_{j\bar k})
 =-\ddc\log\frac{\omega^3}{(\ii)^9\Om\wedge\bar\Om}.
\]
For the smooth positive volume form \(\nu\), we define
\[
 \rho_\nu=-\ddc\log\frac{\nu}{(\ii)^9\Om\wedge\bar\Om}.
\]
Every K\"ahler form satisfying \(\omega^3=\nu\) then has
\(\rho_\omega=\rho_\nu\).

\begin{lemma}\label{lem:positive}
The prescribed volume form \(\nu\) has strictly positive Ricci form
\[
 \rho_\nu=\frac12\ddc\log\bigl((1+t)^2+|u-u^{-1}|^2\bigr)>0.
\]
\end{lemma}
\begin{proof}
We write
\(S=(1+t)^2+|u-u^{-1}|^2=1+\sum_\alpha|f_\alpha|^2\),
where the \(f_\alpha\) are the holomorphic functions
\(\sqrt{2}z_j\), \(z_i z_j\), and \(u-u^{-1}\).
For \(v\in T^{1,0}X\), differentiation and Cauchy--Schwarz give
\[
 (\partial\bar\partial\log S)(v,\bar v)
 =\frac{S\sum_\alpha|\dd f_\alpha(v)|^2
   -\bigl|\sum_\alpha\bar f_\alpha\dd f_\alpha(v)\bigr|^2}{S^2}
 \ge\frac{\sum_\alpha|\dd f_\alpha(v)|^2}{S^2}.
\]
If the left-hand side vanishes, then \(\dd z_j(v)=0\) for every \(j\)
and \(\dd(u-u^{-1})(v)=0\). Differentiating the identity
\(2u=Q(z)+(u-u^{-1})\) on \(X\) gives \(\dd u(v)=0\), so \(v=0\).
\end{proof}

It remains to solve \(\omega^3=\nu\) with \(\omega\) complete.
We use Hein's existence theorem in the formulation of Y.~Li
(\cite{Hein}{}, Proposition~3.1; \cite{Li}{}, Theorem~6.1).
For a background K\"ahler form \(\omega_0\), we write
\[
 f=\log\frac{\nu}{\omega_0^3}.
\]
\begin{theorem}[Hein, in Y.~Li's formulation]\label{thm:Hein}
Let \(\omega_0\) be a K\"ahler form on \(X\). Suppose that, for
some \(0<\alpha<1\), the following conditions hold:
\begin{enumerate}
\item[\textup{(A)}] \((X,g_0)\) is complete and admits a uniform
\(C^{3,\alpha}\) quasi-atlas.
\item[\textup{(B)}] \(f\in C^{2,\alpha}(X)\), and
\( |f(p)|\le C(1+d_{g_0}(p_0,p))^{-\mu}\) for some \(\mu>2\).
\item[\textup{(C)}] \((X,g_0)\) satisfies \(SOB(\beta)\) for some
\(\beta>2\).
\end{enumerate}
Then there are \(0<\bar\alpha\le\alpha\) and a real function
\(\psi\in C^{4,\bar\alpha}(X)\) such that
\[
 \omega=\omega_0+\ddc\psi>0,\qquad \omega^3=e^f\omega_0^3=\nu.
\]
Moreover, \(c\omega_0\le\omega\le C\omega_0\), so \(\omega\)
is complete. If \(f\) is smooth, then \(\psi\) is smooth.
\end{theorem}

The global H\"older norms are measured in the quasi-atlas defined
in Section~\ref{sec:analysis}.
The construction in Appendix~\ref{app:atlas} also verifies the
injectivity-radius and chart-image requirements in Y.~Li's
formulation (\cite{Li}, Definition~5.9).

For Condition C we use the following sufficient form of
\(SOB(\beta)\) (\cite{Li}, Definition~5.7 and the following remark).
We write \(r(p)=d_{g_0}(p_0,p)\). For all sufficiently large \(R\)
and all sufficiently distant \(p\),
we require
\begin{equation}\label{eq:SOBrequirements}
 \Vol B(p_0,R)\le CR^\beta,\qquad
 \Vol B(p,r(p)/2)\ge cr(p)^\beta,\qquad
 \Ric(g_0)\ge-Cr^{-2}g_0.
\end{equation}
We also require that, for some fixed small \(\delta>0\), any two
points with \(r(p)=r(q)=R\) can be joined in
\(\{(1-\delta)R<r<(1+\delta)R\}\).
These conditions give the weighted Sobolev inequalities used in
Hein's theorem. Since \(\beta>2\), no integral normalization of
\(e^f-1\) is required.

Section~\ref{sec:background} constructs the background metric
and proves completeness. Section~\ref{sec:analysis} finishes
the verification of Conditions A--C.

\section{Construction of the background metric}
\label{sec:background}\label{sec:model}\label{sec:distance}

The potentials in Appendix~\ref{app:model} give a background form
\(\omega_0\) whose associated metric is complete and whose volume form satisfies
\[
 \frac{\omega_0^3}{\nu}=1+O((\log t)^{-7/4})
 \qquad(t\to\infty).
\]

\subsection{Outer and inner potentials}\label{subsec:potentials}

On the ends of the base \(\C^*\), we write
\(x=|\log|u|^2|\). The limit \(x\to\infty\) corresponds to
\(|u|\to0\) or \(|u|\to\infty\).

On either end, \(P(u)\ne0\). For fixed \(u\), we choose
\(\sqrt{P(u)}\) and normalize the fibre by \(y=z/\sqrt{P(u)}\).
Then \(Q(y)=1\), and every \(\SO(3)\)-orbit has a representative
\[
 y=\bigl(\cosh(s/2),\ii\sinh(s/2),0\bigr),\qquad s\ge0.
\]
Since \(|y|^2=\cosh s\), this defines the radial parameter
\[
 s=\arcosh\frac{t}{|P(u)|},\qquad t=|P(u)|\cosh s.
\]
The level \(s=0\) is the two-sphere
\[
 \{z=\sqrt{P(u)}\,\eta:\eta\in S^2\subset\R^3\},
\]
where \(t=|P(u)|\) is minimal. For fixed \(u\),
\(s\to\infty\) if and only if \(|z|\to\infty\).
Near \(s=0\),
\[
 \cosh s=1+2|\operatorname{Im}y|^2,\qquad
 s\sim2|\operatorname{Im}y|.
\]
Thus \(s\) is not smooth at \(s=0\), whereas
\(\cosh s=t/|P(u)|\) is smooth.

\medskip
\noindent\phantomsection\label{item:outer-potential}\textbf{Outer potential.}
By Lemma~\ref{lem:F}, \(F=F(\xi,r)\) is smooth
on the cone \(r>|\xi|/2\), even in \(\xi\), and homogeneous
of degree \(5/3\). It satisfies
\begin{equation}\label{eq:model-outer-equation}
 D^2F>0,\qquad F_r>0,\qquad 4F_r\det D^2F=1,
\end{equation}
where \(D^2F\) is the real Hessian in \((\xi,r)\).
We define the outer potential by
\[
 K_o(u,z)=F(\xi,r),
\]
where \(\xi=\log|u|^2\) and \(r=\log(2t)-\arsinh1\).
The ambient complex Hessian has a base--radial block congruent to
\(D^2F\) and, in the two directions orthogonal to \(z\), the block
\((F_r/t)I_2\). By \eqref{eq:model-outer-equation},
\(K_o\) is strictly plurisubharmonic on \(X\) where \(r>|\xi|/2\).
The residue formula expresses the normalized volume density on \(X\)
as the ambient Hessian determinant times the squared norm of the conormal
\(\dd(Q-P)\) in that metric. The calculation in
Appendix~\ref{app:gluing-volume} gives, for fixed \(u\) as \(t\to\infty\),
\[
 \frac{(\ddc K_o)^3}{6(\ii)^9\Om\wedge\bar\Om}
 \sim \frac{F_r^2\det D^2F}{t^3}\,\frac{4t^2}{F_r}
 =\frac{4F_r\det D^2F}{t}
 =\frac1t.
\]
Since \(\nu/(6(\ii)^9\Om\wedge\bar\Om)\sim 1/t\) for fixed \(u\),
this gives \((\ddc K_o)^3\sim\nu\).
See Lemma~\ref{lem:volumeerror} for the uniform error estimate.

\smallskip
Near the cone boundary \(r=|\xi|/2\), we use the expansion
\eqref{eq:Fexp} from Appendix~\ref{app:model}. For
\(x=|\xi|>0\), we set \(s_0=r-x/2\). As \(s_0/x\to0^+\),
\begin{equation}\label{eq:model-outer-expansion}
 \begin{gathered}
 F(\xi,r)=H_0x^{5/3}+B_0C_\Phi x^{1/6}s_0^{3/2}
   +\frac{B_0^2C_\Phi^2}{50H_0}x^{-4/3}s_0^3
   +O(x^{-17/6}s_0^{9/2}),\\
 B_0=\frac3{\sqrt{40H_0}},\quad C_\Phi=\frac{2\sqrt2}{3},
 \end{gathered}
\end{equation}
where \(H_0>0\) is given by \eqref{eq:constants}.
On either base end, \(t=|P(u)|\cosh s\) and
\(|P(u)|=e^{x/2}(1+O(e^{-x}))\) give
\[
 s_0=\log(2\cosh s)-\arsinh1+O(e^{-x}),
\]
uniformly in \(\arg u\). Since \(\log2-\arsinh1<0\), the
outer domain excludes a neighbourhood of the minimal two-sphere
\(s=0\) for large \(x\). For large \(s\),
\[
 s_0=s-\arsinh1+O(e^{-x})+O(e^{-2s}).
\]

\medskip
\noindent\phantomsection\label{item:inner-potential}\textbf{Inner potential.}
The radial functions constructed in Appendix~\ref{app:model},
\(\Phi(s)\) and \(\Psi(s)\), are smooth in \(\cosh s\) for
\(s\ge0\). The form \(\ddc\Phi\) defines a K\"ahler metric \(g_Q\)
on the quadric \(Q(y)=1\). We define the inner potential by
\begin{equation}\label{eq:model-inner-potential}
 K_i(u,z)=H_0x^{5/3}+B_0x^{1/6}\Phi(s)+B_0x^{-4/3}\Psi(s),
\end{equation}
where \(x=|\log|u|^2|\) and \(s=\arcosh(t/|P(u)|)\).

\smallskip
Lemma~\ref{lem:inner} implies that \(K_i\) is smooth and
strictly plurisubharmonic on \(0\le s\le3x^{1/4}\) for
sufficiently large \(x\), including at \(s=0\).
Its Riemannian metric is uniformly comparable to
\begin{equation}\label{eq:model-metric}
 g_*=\frac{20H_0}{9}x^{-1/3}|\dd\log u|^2+B_0x^{1/6}g_Q,
\end{equation}
where, in the local product coordinates \((u,y)\) with \(y=z/\sqrt{P(u)}\),
\(|\dd\log u|^2=|\dd u|^2/|u|^2\) is the cylindrical metric on the
base.
\smallskip
For the model volume form
\[
 \nu_i=\frac{6(\ii)^9\Om\wedge\bar\Om}{\sqrt{t^2+|P(u)|^2}},
\]
the same lemma gives, uniformly on \(0\le s\le3x^{1/4}\),
\begin{equation}\label{eq:model-volume-improvement}
 \frac{(\ddc K_i)^3}{\nu_i}
 =1+O\!\left(\frac{(1+s)^3}{x^3}\right).
\end{equation}
For fixed \(s\), the relative volume error of the uncorrected
inner potential contains a generally nonzero \(x^{-3/2}\) term
(see \eqref{eq:uncorrected-error}). Its decay is too slow for
Condition B. We choose the correction \(\Psi\) through
\eqref{eq:model-correction-equation} to cancel it.
\begin{remark}\label{rem:model-domains}
For sufficiently large \(x\), the outer domain \(r>|\xi|/2\)
over either base end can be written as \(s>s_o(u)\), where
\[
 s_o(u)=\arcosh\frac{e^{x/2+\arsinh1}}{2|P(u)|}
 =\arcosh\frac{1+\sqrt2}{2}+O(e^{-x}).
\]
The error is uniform in \(\arg u\). The inner estimates hold on
\(0\le s\le3x^{1/4}\), so both potentials are smooth and
strictly plurisubharmonic on the overlap
\[
 s_o(u)<s\le3x^{1/4}.
\]
These regions are shown in Figure~\ref{fig:model-domains}.
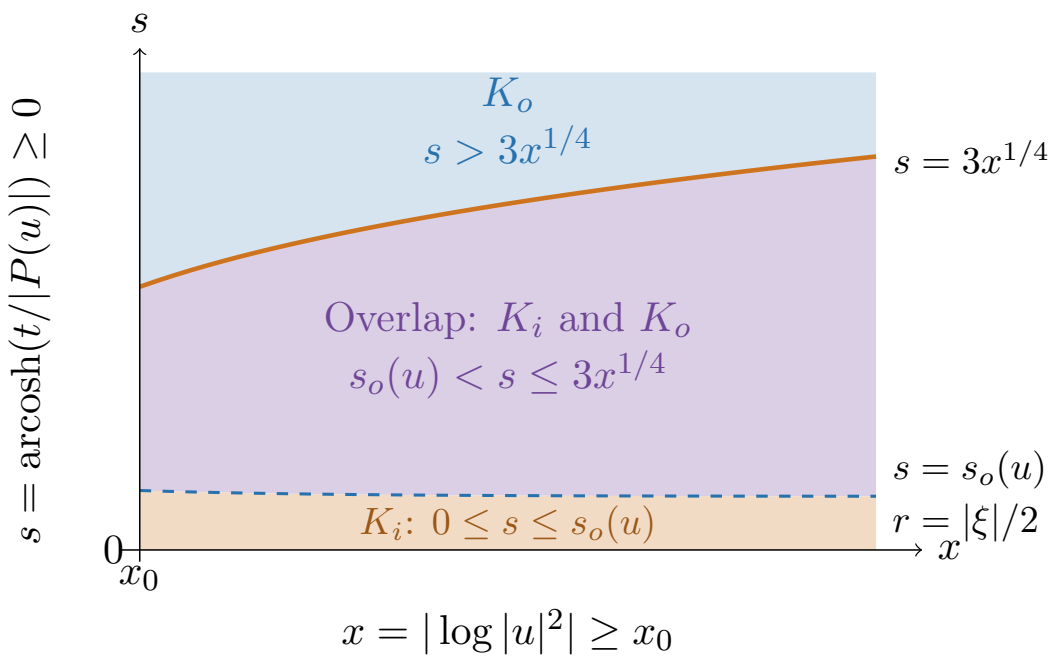
\begin{figure}[H]
\centering
\begin{minipage}[c]{0.48\textwidth}
\centering
\resizebox{\linewidth}{!}{
\begingroup
\definecolor{domainblue}{RGB}{38,113,174}
\definecolor{domainorange}{RGB}{206,116,30}
\definecolor{domainpurple}{RGB}{111,69,151}
\begin{tikzpicture}[x=1.7cm,y=0.9cm,font=\small]
 \fill[domainblue!19]
   (1,4.9)--(5,4.9)
   --plot[domain=5:1,samples=45] (\x,{2.7*pow(\x,0.25)})--cycle;
 \fill[domainpurple!26]
   plot[domain=1:5,samples=45] (\x,{0.55+0.06*exp(1-\x)})
   --plot[domain=5:1,samples=45] (\x,{2.7*pow(\x,0.25)})--cycle;
 \fill[domainorange!26]
   (1,0)--(5,0)
   --plot[domain=5:1,samples=45] (\x,{0.55+0.06*exp(1-\x)})--cycle;
 \draw[domainorange,very thick]
   plot[domain=1:5,samples=45] (\x,{2.7*pow(\x,0.25)})
   node[right,black,font=\footnotesize] {$s=3x^{1/4}$};
 \draw[domainblue,thick,dashed]
   plot[domain=1:5,samples=45] (\x,{0.55+0.06*exp(1-\x)})
   node[right,black,align=left,font=\footnotesize]
   {$s=s_o(u)$\\[2pt]$r=|\xi|/2$};
 \draw[->] (0.9,0)--(5.25,0) node[right] {$x$};
 \draw[->] (1,-0.12)--(1,5.15) node[above] {$s$};
 \node[below] at (1,0) {$x_0$};
 \node[left] at (1,0) {$0$};
 \node[below] at (3,-0.45) {$x=|\log|u|^2|\ge x_0$};
 \node[rotate=90] at (0.4,2.35)
   {$s=\arcosh(t/|P(u)|)\ge0$};
 \node[domainblue,align=center] at (3,4.4)
   {$K_o$\\[2pt]$s>3x^{1/4}$};
 \node[domainpurple,align=center] at (3,2.05)
   {Overlap: $K_i$ and $K_o$\\[2pt]$s_o(u)<s\le3x^{1/4}$};
 \node[domainorange!75!black,font=\footnotesize] at (3,0.26)
   {$K_i$: $0\le s\le s_o(u)$};
\end{tikzpicture}
\endgroup}
\end{minipage}\hspace{0.04\textwidth}%
\begin{minipage}[c]{0.40\textwidth}
\setlength{\captionindent}{0pt}
\setlength{\abovecaptionskip}{0pt}
\raggedright
\caption{Model regions for \(x\ge x_0\).}
\label{fig:model-domains}
\small
\definecolor{domainblue}{RGB}{38,113,174}
\definecolor{domainorange}{RGB}{206,116,30}
\definecolor{domainpurple}{RGB}{111,69,151}
\newcommand{\domainswatch}[1]{%
 \tikz[baseline=-0.5ex]\draw[draw=black!25,fill=#1]
 (0,0) rectangle (0.3,0.2);}
\begin{itemize}
\setlength{\itemsep}{2pt}
\setlength{\parsep}{0pt}
\item \mbox{Inner: \domainswatch{domainorange!26}\,\domainswatch{domainpurple!26}}
\item \mbox{Outer: \domainswatch{domainblue!19}\,\domainswatch{domainpurple!26}}
\end{itemize}
\end{minipage}
\end{figure}

\smallskip
Since \(s_o(u)\) is bounded, for sufficiently large \(x\),
\[
 [x^{1/4},2x^{1/4}]\subset(s_o(u),3x^{1/4}].
\]
On this strip, \(s\asymp x^{1/4}\to\infty\) and
\(s_0/x=O(1/x^{3/4})\to0^+\). Hence the inner asymptotics
\eqref{eq:gexp}, \eqref{eq:psiexp} and the outer boundary
expansion \eqref{eq:model-outer-expansion} apply simultaneously.
\end{remark}

\subsection{Gluing at infinity and extending over the compact core}

We fix \(x_0\) sufficiently large. By Remark~\ref{rem:model-domains},
both potentials are defined and the inner estimates hold on
\(x^{1/4}\le s\le2x^{1/4}\) whenever \(x>x_0\).
Let \(\vartheta:[0,\infty)\to[0,1]\) be smooth, with
\(\vartheta=0\) on \([0,1]\) and \(\vartheta=1\) on
\([2,\infty)\). On each end \(x>x_0\), we define
\begin{equation}\label{eq:gluing}
 K_\infty(u,z)=
 \begin{cases}
 K_i, & 0\le s\le x^{1/4},\\[2pt]
 K_i+\vartheta(s/x^{1/4})(K_o-K_i),
     & x^{1/4}<s<2x^{1/4},\\[2pt]
 K_o, & s\ge2x^{1/4}.
 \end{cases}
\end{equation}
\smallskip
Since \(|P(u)|\le2\cosh(x/2)\), we can choose \(T>0\) such
that, for all \(u\) with \(x<x_0+1\),
\[
 \log(2T)-\arsinh1>\frac{x_0+1}{2},\qquad
 T>|P(u)|\cosh(2x^{1/4}).
\]
On the central region, we define
\[
 K_\infty=K_o\qquad\text{on }\{x<x_0+1,\ t>T\}.
\]
The choice of \(T\) gives \(r>x/2\) there and
\(s>2x^{1/4}\) on the overlap \(x_0<x<x_0+1\), \(t>T\).
The end formula \eqref{eq:gluing} agrees with the central formula
on the overlap, so the two formulas define \(K_\infty\) outside
the compact set \(\{t\le T\}\) (see Lemma~\ref{lem:geometry}).
The central formula also applies at \(P(u)=0\), since \(K_o\)
does not involve \(s\).

\begin{lemma}\label{lem:gluing}
For sufficiently large \(x_0\) and \(T\) as above, the potential
\(K_\infty\) is smooth and strictly plurisubharmonic on
\(\{t>T\}\). Moreover,
\[
 K_\infty\asymp(\log t)^{5/3}
 \qquad\text{uniformly as }t\to\infty.
\]
\end{lemma}
\begin{proof}
The formulas are independent of local choices of \(\log u\)
and \(\sqrt P\). Near \(s=0\),
\(K_\infty=K_i\) is smooth in \(x\) and \(\cosh s\) by
Lemma~\ref{lem:inner}. For \(s>0\), smoothness follows from
\eqref{eq:gluing}, since the cutoff is constant on either side
of the transition interval. On the central region,
\(K_\infty=K_o\) is smooth by Lemma~\ref{lem:F}, including
across \(|u|=1\) and near the fibres \(P(u)=0\).

\smallskip
By Lemmas~\ref{lem:inner} and~\ref{lem:F}, it remains to check
positivity on the transition strip \(x^{1/4}\le s\le2x^{1/4}\).
We compare with \(\ddc K_i\), using the matching expansions
and the cutoff derivative bounds:
\begin{equation}\label{eq:gluing-metric-error}
 \|\ddc K_\infty-\ddc K_i\|_{g_*}
 \le C\left(\frac1{x^{9/4}}+\frac1{x^{7/4}}\right)
 \le\frac{C}{x^{7/4}}.
\end{equation}
The weighted derivative calculation is given in
Appendix~\ref{app:gluing-volume}.
By Lemma~\ref{lem:inner}, the least eigenvalue of \(\ddc K_i\)
relative to \(g_*\) is bounded below by a uniform \(c>0\).
We choose \(x_0\) so large that \(C/x_0^{7/4}<c/2\).
The least eigenvalue of \(\ddc K_\infty\) is then at least \(c/2\)
throughout the transition strip.

\smallskip
The bounds \(\Phi(s)=O((1+s)^{3/2})\) and
\(\Psi(s)=O((1+s)^3)\) imply that, uniformly on
\(0\le s\le2x^{1/4}\),
\[
 \frac{K_i}{H_0x^{5/3}}
 =1+O\!\left(\frac{(1+s)^{3/2}}{x^{3/2}}
             +\frac{(1+s)^3}{x^3}\right)
 =1+O\!\left(\frac1{x^{9/8}}\right).
\]
By \eqref{eq:Fexp}, the same bound holds for \(K_o\) on the
transition strip. Since \eqref{eq:gluing} is a convex combination
there, \(K_\infty\asymp x^{5/3}\) throughout the inner and
transition regions. On these regions,
\[
 \log t=\frac x2+\log\cosh s+O(e^{-x})
       =\frac x2+O(x^{1/4})\asymp x.
\]
Thus \(K_\infty\asymp(\log t)^{5/3}\) there. Elsewhere,
\(K_\infty=K_o\), and Lemma~\ref{lem:F} gives
\[
 K_\infty=F(\xi,r)\asymp r^{5/3}
 \asymp(\log t)^{5/3},\qquad
 r=\log(2t)-\arsinh1.\qedhere
\]
\end{proof}

\medskip
To extend \(K_\infty\) across the compact core, we set
\[
 \chi=\frac12\log\bigl((1+t)^2+|u-u^{-1}|^2\bigr),
\]
so that \(\ddc\chi=\rho_\nu>0\) by Lemma~\ref{lem:positive}.
We choose a relatively compact domain \(U\subset X\) containing
\(\{t\le T\}\), with smooth boundary. The potential \(K_\infty\)
is defined near \(X\setminus U\). By compactness of \(\partial U\),
we choose \(C\) such that
\[
 \chi+C>K_\infty+1\qquad\text{near }\partial U.
\]
The growth estimates \(\chi\asymp\log t\) and
\(K_\infty\asymp(\log t)^{5/3}\) also give
\[
 K_\infty>\chi+C+1\qquad\text{for sufficiently large }t.
\]

Let \(\max_{\mathrm{reg}}\) be a smooth regularized maximum
that agrees with \(\max\) whenever its arguments differ by at
least one. We define
\begin{equation}\label{eq:global-potential}
 K_0(p)=
 \begin{cases}
 \chi(p)+C, & p\in U,\\[2pt]
 \max_{\mathrm{reg}}\!\bigl(\chi(p)+C,K_\infty(p)\bigr),
     & p\in X\setminus U.
 \end{cases}
\end{equation}
Both expressions agree with \(\chi+C\) near \(\partial U\),
which makes \(K_0\) smooth on \(X\). Regularized maxima preserve
strict plurisubharmonicity; hence \(\ddc K_0>0\) on \(X\).
For sufficiently large \(t\),
\(K_0=K_\infty\), so \(K_0\to\infty\) at infinity.
After adding a constant if necessary, we assume \(K_0>0\)
and define the global background form by
\[
 \omega_0=\ddc K_0,\qquad \dd V_0=\frac{\omega_0^3}{6}.
\]
Let \(g_0\) and \(d_{g_0}\) denote the associated Riemannian metric
and distance. Until the solution metric is introduced, all distances
and balls refer to \(g_0\). Outside a compact set,
\(\omega_0=\ddc K_\infty\).

\medskip
\begin{lemma}\label{lem:volumeerror}
The background volume form satisfies
\begin{equation}\label{eq:volumeerror}
 \frac{\omega_0^3}{\nu}
 =1+\wtO{\frac1{(\log t)^{7/4}}}
 \qquad\text{uniformly as }t\to\infty.
\end{equation}
\end{lemma}

On the transition strip \(s\asymp x^{1/4}\), the gluing error
in the volume ratio is \(O(x^{-7/4})\). This suffices for the
pointwise decay in Condition B: \(7/4>5/3\), and the distance is
comparable to \(x^{5/6}\) there by \eqref{eq:distance}.
The weighted remainder is defined and estimated in
Appendix~\ref{app:gluing-volume}.

\subsection{Completeness and distance asymptotics}

We define \(R=\sqrt{4K_0/5}\), a smooth proper function on \(X\).

\begin{proposition}\label{prop:complete}
The background K\"ahler metric \(g_0\) is complete. For every
fixed \(p_0\in X\),
\begin{equation}\label{eq:distance}
 \begin{gathered}
 |\nabla R|_{g_0}=1+o(1),\qquad
 \frac{d_{g_0}(p_0,p)}{R(p)}\longrightarrow1,\\
 d_{g_0}(p_0,p)\asymp(\log t(p))^{5/6},
 \end{gathered}
\end{equation}
uniformly as \(t(p)\to\infty\). Moreover,
\(R|\nabla^2R|_{g_0}\le C\) on \(X\), where the gradient and
Hessian are taken with respect to \(g_0\).
\end{proposition}
\begin{proof}
We use Hermitian norms to compute the model gradients.
In the ambient coordinates \((\zeta,z)\) of
Lemma~\ref{lem:volumeerror}, homogeneity gives
\(D^2F\cdot(\xi,r)^T=\frac23(F_\xi,F_r)^T\). Hence the complex
gradient of \(K_o\) is
\[
 V=\frac32\left(\xi\partial_\zeta+
 r\sum_{j=1}^3 z_j\partial_{z_j}\right),
 \qquad |V|_{g_{\rm amb}}^2=\frac52K_o.
\]
Its orthogonal projection onto \(TX\) is the gradient for the
induced outer metric. Thus
\begin{equation}\label{eq:projection}
 |\partial K_o|_X^2=
 \frac52K_o
 -\frac{9|2rP-\xi(u-u^{-1})|^2}{4N_C},
\end{equation}
where \(N_C\) is the squared conormal norm computed in
Lemma~\ref{lem:volumeerror}. Where \(K_o\) is used, that lemma
gives \(|P|/t=o(1)\) uniformly, so
\(N_C\ge4(t^2-|P|^2)/F_r\ge ct^2/F_r\) for large \(t\).
Using \(|\xi|<2r\), \(F_r\le Cr^{2/3}\), and
\(K_o\asymp r^{5/3}\), we obtain
\[
 0\le1-\frac{2|\partial K_o|_X^2}{5K_o}
 \le\frac{C\log t\,(|P|^2+1)}{t^2}=o(1).
\]
\smallskip
On \(0\le s\le2x^{1/4}\), the formula
\eqref{eq:model-inner-potential} and Lemma~\ref{lem:inner}
give the leading base contribution to \(|\partial K_i|^2\):
\[
 \frac{(\frac53H_0x^{2/3})^2}{\frac{10}{9}H_0x^{-1/3}}
 =\frac52H_0x^{5/3}.
\]
The fibre contribution is
\(O(x^{1/6}(1+s)^{3/2})=o(x^{5/3})\), and the normalized mixed
block tends to zero. Since \(K_i=H_0x^{5/3}(1+o(1))\), we
obtain \(|\partial K_i|^2=\frac52K_i(1+o(1))\) in the inner
metric, uniformly on this region.

\smallskip
On \(x^{1/4}\le s\le2x^{1/4}\),
\eqref{eq:deltader}, \eqref{eq:cutoff-bounds}, and
\eqref{eq:gluing-metric-error} give
\[
 |\partial(K_\infty-K_i)|_{g_*}=o(x^{5/6}),\qquad
 \|\ddc K_\infty-\ddc K_i\|_{g_*}=O(x^{-7/4}).
\]
The first bound controls the gradient change; the second controls
the inverse metric through the uniform comparison with \(g_*\).
The same asymptotic therefore holds for \(K_\infty\).
Since \(K_0-K_\infty\) is constant near infinity, our metric
convention gives
\begin{equation}\label{eq:potential-gradient}
 |\nabla K_0|_{g_0}^2=2|\partial K_0|_{g_0}^2
 =5K_0(1+o(1)).
\end{equation}
Equation~\eqref{eq:potential-gradient} gives
\[
 |\nabla R|_{g_0}^2
 =\frac{|\nabla K_0|_{g_0}^2}{5K_0}=1+o(1).
\]
All remaining norms and covariant derivatives refer to \(g_0\).

\smallskip
The preceding estimate bounds \(|\nabla R|\) near infinity,
and smoothness bounds it on the compact core. Thus
\[
 |R(p)-R(q)|\le C\,d_{g_0}(p,q)\qquad(p,q\in X).
\]
Properness of \(R\) makes every closed \(g_0\)-ball compact,
so \(g_0\) is complete.

\smallskip
We fix \(0<\eps<1\) and choose \(A>R(p_0)\) so large that
\[
 1-\eps\le|\nabla R|\le1+\eps
 \qquad\text{on }\{R\ge A\}.
\]
For \(R(p)>A\), any path from \(p_0\) to \(p\), after its last
exit from \(\{R\le A\}\), has length at least
\((R(p)-A)/(1+\eps)\). For the upper bound, follow
\(-\nabla R/|\nabla R|^2\) from \(p\). The function \(R\)
decreases at unit speed, while the metric speed is at most
\((1-\eps)^{-1}\). The trajectory stays in the compact set
\(\{A\le R\le R(p)\}\) until it reaches \(\{R=A\}\).
Since this level set has bounded distance from \(p_0\),
\[
 \frac{R(p)-A}{1+\eps}
 \le d_{g_0}(p_0,p)\le
 \frac{R(p)-A}{1-\eps}+C_\eps.
\]
After dividing by \(R(p)\), we let \(p\to\infty\) and then
\(\eps\downarrow0\). The distance ratio in \eqref{eq:distance}
follows uniformly at infinity. Lemma~\ref{lem:gluing} gives
\(R\asymp(\log t)^{5/6}\), which proves the distance growth.

\smallskip
Lemma~\ref{lem:radius-hessian} in
Appendix~\ref{app:atlas} gives \(|\nabla^2K_0|\le C\).
Differentiating \(R^2=4K_0/5\) yields
\(R\nabla^2R=\frac25\nabla^2K_0-\dd R\otimes\dd R\), so
\[
 R|\nabla^2R|
 \le\frac25|\nabla^2K_0|+|\nabla R|^2
 \le C.\qedhere
\]
\end{proof}

\medskip
Lemma~\ref{lem:volumeerror} and \eqref{eq:distance} also give
\[
 \frac{\omega_0^3}{\nu}-1
 =O\bigl(d_{g_0}(p_0,p)^{-21/10}\bigr)
 \qquad\text{as }p\to\infty.
\]

\section{Geometric estimates and the Monge--Amp\`ere equation}
\label{sec:analysis}\label{sec:geometry}\label{sec:atlas}
\label{sec:volume}\label{sec:sobolev}\label{sec:MA}

We set \(f=\log(\nu/\omega_0^3)\) and seek a real-valued function
\(\psi\) satisfying
\[
 (\omega_0+\ddc\psi)^3=e^f\omega_0^3=\nu.
\]
We fix \(p_0\in X\). Throughout this section, balls and distances
refer to \(g_0\), and \(\dd V_0=\omega_0^3/6\).

\subsection{Local regularity and Ricci curvature bounds}

\begin{definition}[Quasi-atlas]
Let \((M^n,g)\) be a K\"ahler manifold, and let \(B\subset\C^n\)
be the unit ball. For an integer \(k\ge2\) and \(0<\alpha<1\), a
\emph{\(C^{k,\alpha}\) quasi-atlas} is a family of holomorphic
local diffeomorphisms \(\{\Phi_p:2B\to M\}_{p\in M}\) such that
\(\Phi_p(0)=p\) and
\[
 C^{-1}g_{\mathrm{Euc}}\le\Phi_p^*g\le Cg_{\mathrm{Euc}},\qquad
 \|\Phi_p^*g\|_{C^{k,\alpha}(2B)}\le C_{k,\alpha},
\]
where \(C\ge1\) and \(C_{k,\alpha}>0\) are independent of \(p\).
\end{definition}

For a fixed quasi-atlas, we define
\[
 \|u\|_{C^{\ell,\beta}(M)}
 :=\sup_{p\in M}\|u\circ\Phi_p\|_{C^{\ell,\beta}(B)},
 \qquad \ell\in\Z_{\ge0},\quad 0<\beta<1.
\]
The chart norms are Euclidean; \(C^{\ell,\beta}(M)\) consists
of the functions with finite norm.

\smallskip
\begin{proposition}\label{prop:atlas}
The manifold \((X,g_0)\) admits a quasi-atlas
\(\{\Phi_p\}_{p\in X}\) of class \(C^{k,\alpha}\) for every
integer \(k\ge2\) and every \(0<\alpha<1\). For each integer
\(k\ge0\),
\begin{equation}\label{eq:fchart}
 \|f\circ\Phi_p\|_{C^k(B)}
 \le\frac{C_k}{(1+\log t(p))^{7/4}}
\end{equation}
outside a fixed compact set, with \(C_k\) independent of \(p\).
Consequently, for every \(0<\alpha<1\),
\[
 f\in C^{2,\alpha}(X),\qquad
 |f(p)|\le\frac{C}{\bigl(1+d_{g_0}(p_0,p)\bigr)^{21/10}}
 \quad(p\in X).
\]
\end{proposition}
\begin{proof}
Appendix~\ref{app:atlas} constructs the charts and proves
\eqref{eq:fchart} from Lemma~\ref{lem:volumeerror}. Taking
\(k=3\) gives the \(C^{2,\alpha}\) bound. By
\eqref{eq:distance},
\[
 (\log t(p))^{-7/4}\asymp d_{g_0}(p_0,p)^{-21/10}
 \qquad(p\to\infty),
\]
so \(k=0\) gives the pointwise decay. Both bounds extend to
all of \(X\) by smoothness in the finitely many charts covering
the compact part.
\end{proof}

Appendix~\ref{app:atlas} also checks the injectivity-radius and
image conditions in Y.~Li's formulation of a quasi-atlas
(see \cite{Li}{}, Definition~5.9).

\smallskip
Uniform ellipticity, the \(C^2\) metric bounds, and
\eqref{eq:fchart} with \(k=2\) give
\[
 |\operatorname{Rm}_{g_0}|\le C,\qquad
 |\ddc f|_{g_0}\le
 C\bigl(1+d_{g_0}(p_0,\cdot)\bigr)^{-21/10}.
\]
Since \(\rho_\nu>0\) by Lemma~\ref{lem:positive},
\begin{equation}\label{eq:riclower}
 \rho_{\omega_0}=\rho_\nu+\ddc f
 \ge-C\bigl(1+d_{g_0}(p_0,\cdot)\bigr)^{-21/10}\omega_0
 \ge-C\bigl(1+d_{g_0}(p_0,\cdot)\bigr)^{-2}\omega_0.
\end{equation}

\subsection{Annular connectivity and volume growth}

\begin{lemma}\label{lem:oneend}
The manifold \(X\) has one end, and \(\{R=s\}\) is connected
for all sufficiently large \(s\).
\end{lemma}
\begin{proof}
For \(T>2\), the image of \(\{t=T\}\) under \(\pi_0\) is
\[
 \left\{u\in\C^*:|\xi|\le
 \arcosh\!\left(\frac{T^2}{2}-\cos(2\theta)\right)\right\},
 \qquad \xi=\log|u|^2,\quad\theta=\arg u,
\]
because \(|P(u)|^2=2\cosh\xi+2\cos(2\theta)\). This is a
compact connected annulus. Each fibre of the projection is a
single \(\SO(3)\)-orbit and hence connected. Properness of \(t\)
(Lemma~\ref{lem:geometry}) makes the projection a compact
surjection, so \(\{t=T\}\) is connected.

The map
\(t:\{T_1\le t\le T_2\}\to[T_1,T_2]\) is likewise a compact
surjection with connected fibres when \(T_2>T_1>2\).
Thus these regions, and their increasing union \(\{t>T\}\),
are connected. Since the sublevel sets of \(t\) exhaust \(X\),
the manifold has one end.

\smallskip
By \eqref{eq:distance}, \(R\) is proper and has no critical
points outside a compact set. The vector field
\[
 V=\frac{\nabla R}{|\nabla R|_{g_0}^2},\qquad V(R)=1,
\]
has a flow that identifies \(\{R>A\}\) with
\(\{R=A\}\times(A,\infty)\) for large \(A\).
Each connected component of \(\{R=A\}\) would give an end of
\(X\). Hence this level, and every larger level, is connected.
\end{proof}

\smallskip
For each fixed \(0<\delta<1\), any two points \(p,q\) with
\(d_{g_0}(p_0,p)=d_{g_0}(p_0,q)=r\) can therefore be joined in
\[
 \{(1-\delta)r<d_{g_0}(p_0,\cdot)<(1+\delta)r\}
\]
when \(r\) is sufficiently large. Indeed, choose
\(0<\eps<\delta/4\). Equation~\eqref{eq:distance} gives
\(R(p),R(q)\in[(1-\eps)r,(1+\eps)r]\). Follow the trajectories
of \(V\) from both points to \(\{R=r\}\), and join the resulting
points within that connected level. The path stays in
\[
 \{(1-\eps)r\le R\le(1+\eps)r\}
 \subset
 \{(1-\delta)r<d_{g_0}(p_0,\cdot)<(1+\delta)r\},
\]
where the inclusion follows again from the uniform distance ratio
in \eqref{eq:distance}.

\medskip
\begin{lemma}\label{lem:pushvolume}\label{prop:remote}
There exist constants \(c,C>0\) such that, for all sufficiently
large \(r\),
\begin{equation}\label{eq:volume-growth}
 cr^{12/5}\le\Vol_{g_0}B(p_0,r)\le Cr^{12/5}.
\end{equation}
Moreover, for every \(p\) sufficiently far from \(p_0\),
\begin{equation}\label{eq:remote}
 \Vol_{g_0}B\bigl(p,\tfrac12d_{g_0}(p_0,p)\bigr)
 \ge c\,d_{g_0}(p_0,p)^{12/5}.
\end{equation}
\end{lemma}
\begin{proof}
We first express \(\nu\) in logarithmic and angular coordinates.
For \(P\ne0\) and \(t>|P|\), set \(y=z/\sqrt P\) and use the
\(\SO(3)\)-orbit representative
\[
 y=(\cosh(s/2),\ii\sinh(s/2),0).
\]
The Euclidean metric coefficients in the radial and three orbit
directions are
\[
 \frac{\cosh s}{4},\quad\cosh s,\quad
 \cosh^2(s/2),\quad\sinh^2(s/2).
\]
Their volume density is \(\cosh s\sinh s/4\). Division by
\(|\partial Q|^2=4\cosh s\) gives the residue measure
\(c\sinh s\,\dd s\,\dd\mu_{\SO(3)}\), where
\(\mu_{\SO(3)}\) is a fixed Haar measure.
The product formula \eqref{eq:residueproduct} supplies the factor
\(|P|\,\dd\xi\,\dd\theta\). Since
\(\dd t=|P|\sinh s\,\dd s\) along each fibre, \eqref{eq:nu}
becomes
\begin{equation}\label{eq:pushvolume}
 \nu=c\,\frac{t}{\sqrt{(1+t)^2+|u-u^{-1}|^2}}\,
 \dd\xi\,\dd\log t\,\dd\theta\,\dd\mu_{\SO(3)}.
\end{equation}
The excluded fibres and degenerate orbits have measure zero.
The inequality
\(|u-u^{-1}|^2\le|P|^2+4\le t^2+4\) bounds the density factor
above and below for large \(t\). By \eqref{eq:volumeerror},
\[
 \dd V_0\asymp\nu\asymp
 \dd\xi\,\dd\log t\,\dd\theta\,\dd\mu_{\SO(3)}
 \qquad(t\gg1).
\]

\smallskip
For large \(t\), uniformly in \(\theta\),
\[
 \{|\xi|\le\log t\}\subset\{|P(u)|\le t\}
 \subset\{|\xi|\le2\log t+C\}.
\]
Integrating over these strips and the compact angular variables
gives, for fixed large \(T\) and \(L\to\infty\),
\[
 \nu\{T\le t\le e^L\}
 \asymp\int_{\log T}^L\ell\,\dd\ell\asymp L^2.
\]
Equation~\eqref{eq:distance} gives constants \(c_0,C_0>0\) with
\[
 \{t\le e^{c_0r^{6/5}}\}\subset B(p_0,r)
 \subset\{t\le e^{C_0r^{6/5}}\}
\]
for large \(r\). Together with the preceding volume comparison
and the finite volume of the compact part, this proves
\eqref{eq:volume-growth}.

\medskip
To prove \eqref{eq:remote}, we fix \(p\) far from \(p_0\) and write
\(r=d_{g_0}(p_0,p)\). We construct a square of side
\(\eps\log t(p)\) in the \((\xi,\log t)\)-plane whose full
inverse image in \(X\) lies in \(B(p,r/2)\).

If \(p\) lies on an end with \(s<3x^{1/4}\), first increase
\(s\) to \(3x^{1/4}\) at fixed \(u\). Lemma~\ref{lem:inner},
\eqref{eq:gluing-metric-error}, and the outer expansion
\eqref{eq:Fexp} bound the length by
\[
 Cx^{1/12}\int_0^{3x^{1/4}}\frac{\dd s}{(1+s)^{1/4}}
 \le Cx^{13/48}=o(r).
\]
Here \(r\asymp x^{5/6}\), and the change in \(\log t\) is
\(O(x^{1/4})=o(\log t(p))\). The endpoint lies where
\(\omega_0=\ddc K_o\). This identity already holds in the
central region and when \(s\ge3x^{1/4}\).

\smallskip
At fixed \((\xi,t)\), all base phases are allowed and their
lifts remain in the outer region. Indeed, on either end,
\(2\sinh(x/2)\le |P(u)|\le2\cosh(x/2)\). Since the
preceding step ensures \(s\ge3x^{1/4}\), we have
\[
 \frac{t}{2\cosh(x/2)}
 \ge\tanh(x/2)\cosh(3x^{1/4})
 >\cosh(2x^{1/4})
\]
for large \(x\). Thus every base phase satisfies
\(s>2x^{1/4}\). For bounded \(x\), the same outer-metric
identity holds for large \(t\), by the central prescription
or because \(s>2x^{1/4}\) on the ends. The \(\SO(3)\)
rotations preserve \(u\) and \(t\), so they also stay in
the outer region. To lift a change in \(\theta\), use
\eqref{eq:central-chart} with real \(\lambda\) and
\[
 v=0,\qquad \zeta=\xi/2+\ii\theta,\qquad
 e^{2\lambda}=t+\sqrt{t^2-|P|^2}.
\]
Since \(2\cosh(x/2)/t=o(1)\), these paths and the \(\SO(3)\) rotation
paths have bounded speeds in \((\zeta,z/\sqrt t)\). For each
such path \(\gamma\), the outer metric, \eqref{eq:Fexp}, and
homogeneity on interior subcones give
\[
 |\dot\gamma|_{g_0}^2
 \le C(F_{\xi\xi}+F_{rr}+F_r)\le C(\log t)^{2/3}.
\]
The parameter intervals have bounded length, so these paths have
total length at most \(C(\log t)^{1/3}=o(r)\).
These paths reach the real surface
\[
 \Sigma:\quad u=e^{\xi/2}>0,\qquad
 z_1=\sqrt{(t+P)/2},\quad
 z_2=\ii\sqrt{(t-P)/2},\quad z_3=0.
\]
On \(\Sigma\), we write \(h=\tfrac12D^2F\), with the Hessian
evaluated at \((\xi,\log(2t)-\arsinh1)\). Differentiation gives
\[
 g_0|_\Sigma=h
 +\frac{F_r}{2(t^2-P^2)}(\dd P-P\,\dd\log t)^2.
\]
Since \(P=2\cosh(\xi/2)\), we have \(|P_\xi|\le P/2\).
Cauchy--Schwarz and the inverse-Hessian bound
\(\|F_r(D^2F)^{-1}\|\le C\log t\), proved in the derivation
of \eqref{eq:outer-volume-error}, yield the quadratic-form estimate
\[
 0\le g_0|_\Sigma-h
 \le C(\log t)\frac{P^2}{t^2-P^2}\,h.
\]
On the ends, \(P/t=1/\cosh s\) and \(s\ge2x^{1/4}\)
imply \(\log t\le C(1+s^4)\), so the relative error is
bounded by \(C(1+s^4)/\sinh^2s=o(1)\). On a fixed
interior subcone, \(P/t\le Ct^{-\delta/2}\) for some
\(\delta>0\), giving an error \(O((\log t)t^{-\delta})=o(1)\).
Thus \(g_0|_\Sigma=(1+o(1))h\) uniformly on this region.

\smallskip
By Lemma~\ref{lem:hessiancone},
\[
 \frac12D^2F=(\dd R_F)^2+R_F^2(\dd\sigma)^2,\qquad
 R_F=\sqrt{\frac{4F}{5}},\qquad 0<\sigma<\frac{5\pi}{6}.
\]
The preceding paths have length \(o(r)\), so
\(R_F\asymp(\log t(p))^{5/6}\asymp r\). Fix a small
\(\eta>0\) and move \(\sigma\) to the nearest point of
\([\eta,5\pi/6-\eta]\), keeping \(R_F\) fixed. Along this
arc \(\gamma\), we have \(F=\tfrac54R_F^2=\mathrm{const}\).
The arc stays in the outer region. Evenness and convexity give
\(F_\xi\ge0\) on the positive end. Moving inward decreases
\(x\), while
\[
 \left.\frac{\dd\log t}{\dd x}\right|_F
 =-\frac{F_\xi}{F_r}\le0;
\]
thus \(t\) increases, while \(P=2\cosh(x/2)\) decreases.
Consequently \(s=\arcosh(t/P)\) increases and the threshold
\(2x^{1/4}\) decreases, so the arc stays in the outer branch
wherever the end prescription applies. The central prescription
covers any remaining part of the arc. The negative end is symmetric.
Since \(F\asymp(\log t)^{5/3}\), \(\log t\) remains comparable
to \(\log t(p)\) along the arc. Together with the decrease of
\(P/t\), this preserves the uniform relative error above. Hence
\[
 L_{g_0}(\gamma)\le(1+o(1))R_F\eta\le C\eta r.
\]
We choose \(\eta\), independently of \(p\), so that the total
length from \(p\) to the endpoint \(p_*\in\Sigma\) is less
than \(r/4\) for large \(r\). We also have
\(\log t(p_*)\asymp\log t(p)\).

We take the square centred at \((\xi(p_*),\log t(p_*))\), with side
length \(\eps\log t(p)\), shown in Figure~\ref{fig:remote-volume}.
For \(\eps>0\) sufficiently small relative to \(\eta\), the square
lies in a fixed interior subcone and \(\log t\asymp\log t(p)\)
throughout it.

\begin{figure}[ht]
\centering
\resizebox{0.52\textwidth}{!}{
\begingroup
\definecolor{volblue}{RGB}{38,113,174}
\definecolor{volorange}{RGB}{190,100,25}
\begin{tikzpicture}[x=1.05cm,y=0.95cm,font=\small]
  \fill[volblue!7]
    (-4.5,2.55)--(0,0.3)--(4.5,2.55)--(4.5,4.6)--(-4.5,4.6)--cycle;
  \draw[gray,dashed] (-4.5,2.55)--(0,0.3)--(4.5,2.55);
  \draw[->] (-4.7,0)--(4.8,0) node[right] {$\xi$};
  \draw[->] (0,-0.15)--(0,4.85) node[above] {$\log t$};
  \node[gray,rotate=-27] at (-3.2,1.35) {cone boundary};

  \filldraw[fill=volblue!22,draw=volblue,thick]
    (1.0,2.8) rectangle (2.2,4.0);
  \fill[volblue] (1.6,3.4) circle (1.6pt);
  \node[volblue,below left,font=\footnotesize] at (1.6,3.4) {$p_*$};
  \fill (4.0,2.38) circle (2pt);
  \node[anchor=west] at (4.15,2.38) {$p$};
  \draw[gray,thick,dashed,->] (4.0,2.43)--(4.0,2.95);
  \node[anchor=east,font=\footnotesize] at (3.85,2.68) {$o(r)$};
  \draw[volorange,thick,->]
    (4.0,3.0) .. controls (3.8,3.35) and (2.8,3.38) .. (1.65,3.4);
  \fill[volorange] (4.0,3.0) circle (1.8pt);
  \node[volorange,above,font=\footnotesize] at (3.3,3.38) {$F=\mathrm{const}$};

  \draw[volblue,<->] (1.0,4.2)--(2.2,4.2);
  \node[volblue,above,font=\footnotesize] at (1.6,4.2)
    {$\varepsilon\log t(p)$};
  \node[volblue,align=center,font=\footnotesize] at (1.6,2.13)
    {$\mathrm{area}=\varepsilon^2(\log t(p))^2$};
\end{tikzpicture}
\endgroup}
\caption{Projection to \((\xi,\log t)\). The square's inverse
image in \(X\) lies in \(B(p,r/2)\).}
\label{fig:remote-volume}
\end{figure}
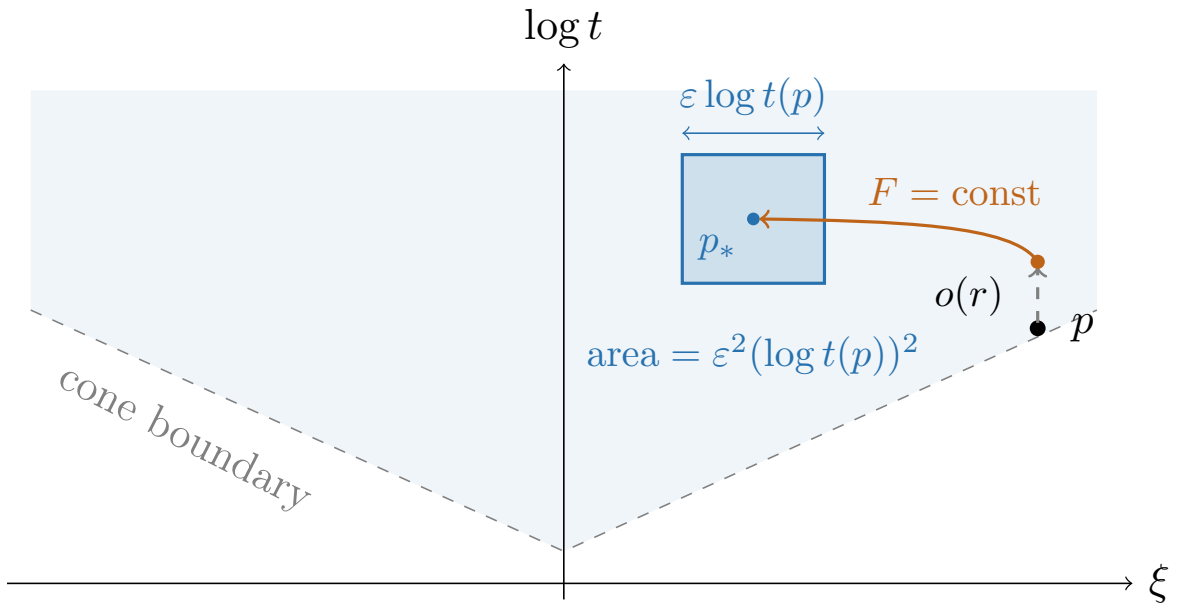

Homogeneity of \(D^2F\) gives
\[
 g_0|_\Sigma\le
 \frac{C_\eta}{(\log t(p))^{1/3}}
 \bigl((\dd\xi)^2+(\dd\log t)^2\bigr).
\]
A segment from the centre to any point of the square therefore
lifts to a path in \(\Sigma\) of length at most
\[
 C_\eta\eps(\log t(p))^{5/6}\le C'_\eta\eps r<r/8,
\]
after fixing \(\eps\) sufficiently small, independently of \(p\).

\smallskip
For some fixed \(0<\delta<2\), the square satisfies
\(|\xi|\le(2-\delta)\log t\). Every base phase is admissible,
since
\[
 |P(u)|\le e^{\xi/2}+e^{-\xi/2}
 \le2t^{1-\delta/2}<t.
\]
For any \(q\in X\) whose \((\xi,\log t)\)-coordinates lie in
the square, join \(p\) to \(p_*\),
follow the segment in \(\Sigma\) to the coordinates
\((\xi(q),\log t(q))\), and then vary the angular variables.
The three length bounds give
\[
 d_{g_0}(p,q)<\frac r4+\frac r8+C(\log t(p))^{1/3}<\frac r2
\]
for large \(r\). Integrating \eqref{eq:pushvolume} over the
square and all angular variables gives
\[
 \Vol_{g_0}B(p,r/2)
 \ge c\eps^2(\log t(p))^2\ge c'r^{12/5},
\]
which proves \eqref{eq:remote}.
\end{proof}

\subsection{Solving the Monge--Amp\`ere equation}

Propositions~\ref{prop:complete} and~\ref{prop:atlas} give
Condition A of Theorem~\ref{thm:Hein}. Proposition~\ref{prop:atlas}
also gives Condition B with \(\mu=21/10>2\).

\smallskip
The annular connectivity proved after Lemma~\ref{lem:oneend}
is the connectivity condition in \(SOB(\beta)\)
(see \cite{Li}{}, Definition~5.7 and the following remark).
For \(r=d_{g_0}(p_0,p)\) sufficiently large, we choose a common
constant \(C\ge2\). Equations~\eqref{eq:volume-growth},
\eqref{eq:remote}, and \eqref{eq:riclower} give
\[
 \begin{aligned}
 \Vol_{g_0}B(p_0,r)&\le Cr^{12/5},\\
 \Vol_{g_0}B\bigl(p,(1-C^{-1})r\bigr)
 &\ge\Vol_{g_0}B(p,r/2)\ge C^{-1}r^{12/5},\\
 \Ric_{g_0}(p)&\ge-Cr^{-2}g_0(p).
 \end{aligned}
\]
Thus \((X,g_0)\) satisfies Condition C with \(\beta=12/5>2\).

\medskip
\begin{proof}[Proof of Theorem~\ref{thm:main}]
Theorem~\ref{thm:Hein} gives a real-valued
\(\psi\in C^{4,\bar\alpha}(X)\), where
\(0<\bar\alpha\le\alpha\), satisfying
\[
 \omega=\omega_0+\ddc\psi>0,\qquad
 \omega^3=e^f\omega_0^3=\nu.
\]

\smallskip
The global chart bounds imply \(|\ddc\psi|_{g_0}\le C\).
If \(\lambda_1,\lambda_2,\lambda_3\) are the eigenvalues of
\(\omega\) relative to \(\omega_0\), then
\[
 0<\lambda_i\le C,\qquad
 \prod_{j=1}^3\lambda_j=e^f,\qquad
 \lambda_i\ge C^{-2}e^{-\|f\|_{L^\infty(X)}}.
\]
Hence
\begin{equation}\label{eq:comparison}
 c\omega_0\le\omega\le C\omega_0.
\end{equation}
Completeness of \(g_0\) and the lower bound imply completeness
of \(\omega\). The smoothness of \(f\) and \(\omega_0\), together
with uniform ellipticity and local regularity for the Monge--Amp\`ere
equation, gives \(\psi\in C^\infty(X)\). The volume equation
\(\omega^3=\nu\) and Lemma~\ref{lem:positive} give
\[
 \rho_\omega=\rho_\nu
 =\frac12\ddc\log\bigl((1+t)^2+|u-u^{-1}|^2\bigr)>0.
\]

\smallskip
The projection \(\pi_0:X\to\C^*\) has the holomorphic section
\begin{equation}\label{eq:section}
 \sigma(u)=\left(u,\frac{P(u)+1}{2},\frac{P(u)-1}{2\ii},0\right),
 \qquad \pi_0\circ\sigma=\id.
\end{equation}
Thus
\((\pi_0)_*\circ\sigma_*=\id\) on homology and fundamental
groups, giving
\[
 \QQ\cong H_1(\C^*;\QQ)\lhook\joinrel\longrightarrow H_1(X;\QQ),
 \qquad
 (\pi_0)_*:\pi_1(X)\twoheadrightarrow\pi_1(\C^*)\cong\Z.
\]
\end{proof}

\section{Noncompactifiability and divergent Ricci integrals}
\label{sec:topology}

Consider the map in Corollary~\ref{thm:cover},
\[
 \pi:M=\{(Z,w)\in\C^3\times\C:Q(Z)=\sin w\}\longrightarrow X,
 \qquad \pi(Z,w)=(-\ii e^{\ii w},\sqrt2Z).
\]

\begin{lemma}\label{lem:covergeometry}\label{lem:cycles}
The manifold \(M\) is smooth and connected, and \(\pi\) is a
regular holomorphic covering with deck group
\(\langle T\rangle\cong\Z\), where \(T(Z,w)=(Z,w+2\pi)\).
Moreover,
\[
 \dim_\QQ H_3(M;\QQ)=\infty.
\]
\end{lemma}
\begin{proof}
The identities
\[
 P(-\ii e^{\ii w})=2\sin w,\qquad Q(\sqrt2Z)=2Q(Z)
\]
identify \(\pi\) with the pullback of the universal covering
\(w\mapsto-\ii e^{\ii w}\) of \(\C^*\) along \(\pi_0:X\to\C^*\).
The section in \eqref{eq:section} makes
\((\pi_0)_*:\pi_1(X)\to\pi_1(\C^*)\) surjective, so the
pullback is connected and its deck group is generated by \(T\).
Since \(X\) is smooth and \(\pi\) is locally biholomorphic,
\(M\) is smooth.

\smallskip
For \(k\in\Z\), we define
\[
 S_k=\left\{(\ii^k\sqrt{|\sin w|}\,\eta,w):
 \eta\in S^2\subset\R^3,\quad k\pi\le w\le(k+1)\pi\right\}.
\]
On this interval, \(\sin w=(-1)^k|\sin w|\), so \(S_k\subset M\).
Over the open interval, \(S_k\) is a product with \(S^2\);
at each endpoint the sphere collapses to a point. For
\(v\in\R^3\) sufficiently small, the endpoint charts are
\[
 v\longmapsto(\ii^kv,k\pi+\arcsin|v|^2),\qquad
 v\longmapsto(\ii^kv,(k+1)\pi-\arcsin|v|^2).
\]
These charts show that \(S_k\) is a smooth embedded three-sphere.
Its projection to the \(w\)-plane is the interval
\([k\pi,(k+1)\pi]\), as shown in
Figure~\ref{fig:vanishing-spheres}.
\begin{figure}[H]
 \centering
\begin{tikzpicture}[x=1.2cm,y=0.65cm,font=\small]
 \draw[->,gray] (-0.35,0)--(5.6,0) node[right] {$\operatorname{Re}w$};
 \draw[thick,blue!65!black] (0,0)--(1,0);
 \draw[thick,orange!80!black] (1,0)--(2,0);
 \draw[thick,gray] (2,0)--(3,0);
 \draw[thick,blue!65!black] (3,0)--(4,0);
 \draw[thick,orange!80!black] (4,0)--(5,0);
 \foreach \x/\lab in {0/0,1/\pi,2/2\pi,3/3\pi,4/4\pi,5/5\pi}{
  \fill (\x,0) circle (1.6pt);
  \node[below] at (\x,-0.05) {$\lab$};
 }
 \foreach \x/\lab in {0.5/0,1.5/1,2.5/2,3.5/3,4.5/4}{
  \node[above] at (\x,0.1) {$S_{\lab}$};
 }
 \node at (-0.55,0) {$\cdots$};
 \node[above] at (5.35,0.1) {$\cdots$};
\end{tikzpicture}
 \caption{Projection of the spheres \(S_k\) to the real \(w\)-axis.}
 \label{fig:vanishing-spheres}
\end{figure}
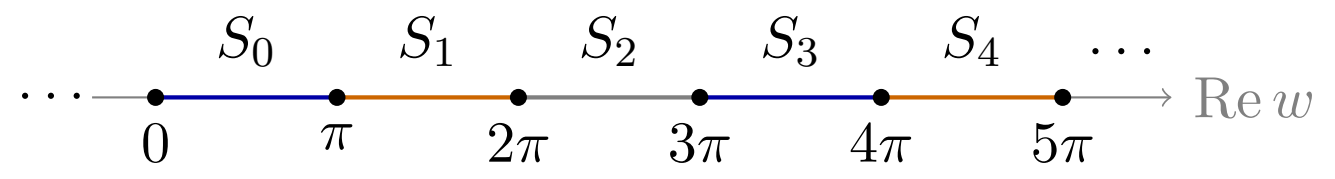

Adjacent spheres meet only at \((0,(k+1)\pi)\). Since
\(\cos((k+1)\pi)\ne0\), \(Z\) gives holomorphic coordinates
on \(M\) near this point, and
\[
 T_{(0,(k+1)\pi)}S_k=\ii^k\R^3,\qquad
 T_{(0,(k+1)\pi)}S_{k+1}=\ii^{k+1}\R^3.
\]
Their intersection is transverse, with intersection number
\(\pm1\), and nonadjacent spheres are disjoint. We select the spheres
\(S_{3j}\) for \(j\ge0\) to obtain independent classes. Among the
test spheres \(S_{3l+1}\), \(l\ge0\), only \(S_{3j+1}\) meets \(S_{3j}\).
We choose orientations so that
\[
 [S_{3j}]\cdot[S_{3l+1}]=\delta_{jl},\qquad j,l\ge0.
\]
Pairing with the Poincar\'e duals of \(S_{3l+1}\), which lie in
\(H_c^3(M;\QQ)\), proves linear independence:
\[
 \sum_{j=0}^N a_j[S_{3j}]=0
 \quad\Longrightarrow\quad
 0=\sum_{j=0}^N a_j\bigl([S_{3j}]\cdot[S_{3l+1}]\bigr)=a_l
 \quad(0\le l\le N).\qedhere
\]
\end{proof}

\medskip
\begin{proof}[Proof of Corollary~\ref{thm:cover}]
Since \(\pi\) is a holomorphic covering, \(\omega_M=\pi^*\omega\)
is complete and has positive Ricci form
\(\rho_{\omega_M}=\pi^*\rho_\omega\).
A sheet over an evenly covered ball has infinitely many disjoint
deck translates. By invariance under deck transformations, each
has the same positive integral of \(\rho_{\omega_M}^3\). Summing
over these translates gives \(\int_M\rho_{\omega_M}^3=\infty\).

\smallskip
Suppose that \(M\) is an analytic Zariski open subset of a
compact complex manifold. By subanalytic triangulation of the
compactification and its analytic complement (see
\cite{HironakaTriangulation}{}, p.~180), \(M\) has the homotopy type
of a finite CW complex. This contradicts
\(\dim_\QQ H_3(M;\QQ)=\infty\) from Lemma~\ref{lem:cycles}.
\end{proof}

\medskip
\begin{proof}[Proof of Corollary~\ref{cor:ricci-integrals}]
We take \(N_{n,k}=M\times\mathbb P^{k-3}\times\C^{n-k}\) with the
product metric, using \(\omega_M\) from Corollary~\ref{thm:cover},
the Fubini--Study metric, and the planar metric
\[
 g_\C=\frac{|\dd\zeta|^2}{\sqrt{1+|\zeta|^2}}.
\]
The planar metric has positive Ricci form
\(\tfrac12\ddc\log(1+|\zeta|^2)\). For large Euclidean radius \(T\),
the radial distance is comparable to \(T^{1/2}\) and the disc area
to \(T\), by integrating the radial length and area elements.
Thus \(g_\C\) is complete and its balls of radius \(R\) have area
comparable to \(R^2\) at any fixed centre. The product is complete
and has positive Ricci curvature.

All terms in the wedge expansion are nonnegative. Retain three
Ricci factors from \(M\), \(k-3\) from projective space, and one
K\"ahler factor from each planar factor. By the product distance
formula, for a fixed centre \(p\) and sufficiently small \(c>0\) the ball
\(B_{\omega_{n,k}}(p,R)\) contains the product of the \(M\)-ball
\(B_{\omega_M}(p_M,cR)\), the whole projective factor, and the planar
balls of radius \(cR\), for all large \(R\). Integrating the retained
term over this product gives
\[
 R^{2k-2n}\int_{B_{\omega_{n,k}}(p,R)}
       \rho_{\omega_{n,k}}^k\wedge\omega_{n,k}^{n-k}
 \ge C^{-1}\int_{B_{\omega_M}(p_M,cR)}\rho_{\omega_M}^3
 \longrightarrow\infty.
\]
Here the projective factor contributes a positive constant, the
planar factors contribute a factor comparable to \(R^{2(n-k)}\),
and a zero-dimensional factor contributes \(1\).

Fixing the other coordinates defines an inclusion of \(M\).
Projection onto \(M\) is a left inverse, so \(H_3(N_{n,k};\QQ)\)
contains a copy of \(H_3(M;\QQ)\) and is infinite-dimensional.
The finite CW type obstruction used for \(M\) again excludes
an analytic compactification.
\end{proof}

\appendix
\section{Construction and estimates for the model potentials}\label{app:model}

We construct the functions \(F\), \(\Phi\), and \(\Psi\) used in
Subsection~\ref{subsec:potentials}. Throughout the appendices,
\[
 a=\frac53,\qquad p=\frac32,\qquad d=a-p=\frac16.
\]
Weighted remainders use the conventions in Appendix~\ref{app:gluing-volume}.

\subsection{The outer potential and its boundary expansion}

We solve \eqref{eq:model-outer-equation} on \(r>|\xi|/2\)
with \(F\) even in \(\xi\) and homogeneous. If \(F\) has degree
\(a\), then \(F_r\det D^2F\) has degree \(3a-5\), so the
equation requires \(a=5/3\).
Let \(B(\cdot,\cdot)\) denote the Euler beta function. We set
\begin{equation}\label{eq:constants}
 H_0=\frac{9}{5(16/9)^{1/3}B(1/2,2/3)}.
\end{equation}
We define
\begin{equation}\label{eq:I-def}
 I(T)=\int_1^T(\sigma^{6/5}-1)^{-1/3}\dd\sigma,\qquad T\ge1.
\end{equation}
The integrand is positive, locally integrable at \(1\), and
asymptotic to \(\sigma^{-2/5}\) at infinity. Thus \(I\) is continuous
and strictly increasing
from \([1,\infty)\) onto \([0,\infty)\). For each \(y\ge0\), let
\(T\ge1\) be the unique solution of
\begin{equation}\label{eq:varphi-def}
 y=(16/9)^{1/3}H_0 I(T),\qquad \varphi(y)=H_0T.
\end{equation}
Then \(\varphi\) is continuous on
\([0,\infty)\), satisfies \(\varphi(0)=H_0\), and is smooth on
\((0,\infty)\). Implicit differentiation gives
\begin{equation}\label{eq:varphi-ode}
 (\varphi')^3=\frac9{16}\left[
       \left(\frac{\varphi}{H_0}\right)^{6/5}-1\right],
 \qquad
 \varphi''=\frac9{40H_0\varphi'}
       \left(\frac{\varphi}{H_0}\right)^{1/5}>0.
\end{equation}
In particular, \(\varphi'>0\), and the identities in
\eqref{eq:varphi-ode} yield
\begin{equation}\label{eq:varphi-identity}
 a\varphi\varphi''-(a-1)(\varphi')^2
 =\frac3{8\varphi'}.
\end{equation}

\smallskip
For \(\xi\ne0\), we set \(x=|\xi|\) and define
\begin{equation}\label{eq:F}
 F(\xi,r)=x^a\varphi\left(\frac rx-\frac12\right),
 \qquad r>x/2.
\end{equation}

\begin{lemma}\label{lem:F}
The function \(F\) defined by \eqref{eq:F} extends to a smooth even
function on \(\{r>|\xi|/2\}\), and this extension is real analytic
across \(\xi=0\). On this domain,
\begin{equation}\label{eq:null}
 \begin{gathered}
 D^2F>0,\qquad F_r>0,\qquad 4F_r\det D^2F=1,\\
 F(\lambda\xi,\lambda r)=\lambda^aF(\xi,r)\quad(\lambda>0).
 \end{gathered}
\end{equation}

Moreover, \(F\asymp r^a\) uniformly on \(\{r>|\xi|/2\}\).
\end{lemma}
\begin{proof}
Evenness and homogeneity follow from \eqref{eq:F}.
On \(\xi>0\), we set \(x=\xi\) and
\[
 y=\frac rx-\frac12,\qquad v=y+\frac12=\frac rx,
 \qquad y_r=\frac1x,\qquad y_x=-\frac vx.
\]
The first derivatives of \(F=x^a\varphi(y)\) are
\begin{equation}\label{eq:F-first-derivatives}
 F_r=x^{a-1}\varphi',\qquad
 F_x=x^{a-1}\bigl(a\varphi-v\varphi'\bigr).
\end{equation}

Differentiating once more gives
\begin{equation}\label{eq:F-second-derivatives}
\begin{aligned}
 F_{rr}&=x^{a-2}\varphi'',\\
 F_{xr}&=x^{a-2}\bigl((a-1)\varphi'-v\varphi''\bigr),\\
 F_{xx}&=x^{a-2}\bigl(a(a-1)\varphi
       -2(a-1)v\varphi'+v^2\varphi''\bigr).
\end{aligned}
\end{equation}
\smallskip
The terms containing \(v\) cancel in the determinant:
\begin{equation}\label{eq:F-determinant}
 \det D^2F
 =(a-1)x^{2a-4}
   \bigl(a\varphi\varphi''-(a-1)(\varphi')^2\bigr).
\end{equation}
Since \(3a-5=0\) and \(a-1=2/3\), equations
\eqref{eq:F-first-derivatives}, \eqref{eq:F-determinant}, and
\eqref{eq:varphi-identity} give
\[
 4F_r\det D^2F
 =\frac32(a-1)x^{3a-5}=1.
\]
By \eqref{eq:varphi-ode}, \(\varphi',\varphi''>0\), so
\(F_r,F_{rr}>0\). The equation gives \(\det D^2F>0\), and hence
\(D^2F>0\). Evenness gives the same conclusions on \(\xi<0\).

\medskip
To extend \(F\) across \(\xi=0\), we fix \(r>0\) and let
\(y\to\infty\). In \eqref{eq:I-def}, we substitute
\(\sigma^{6/5}=\tau^{-1}\) and integrate by parts to obtain
\begin{equation}\label{eq:I-by-parts}
 I(T)=\frac53T^{3/5}(1-T^{-6/5})^{2/3}
 -\frac5{18}\int_{T^{-6/5}}^1
       \tau^{-1/2}(1-\tau)^{-1/3}\dd\tau.
\end{equation}
For sufficiently large \(T\), convergent binomial expansions give
\begin{equation}\label{eq:I-large}
 I(T)=aT^{1/a}-\frac5{18}B(1/2,2/3)+
       \sum_{m\ge1}c_mT^{-(2m-1)/a}.
\end{equation}
The choice \eqref{eq:constants} gives
\[
 (16/9)^{1/3}H_0\frac5{18}B(1/2,2/3)=\frac12.
\]
Thus the constant term cancels in \(r/x=y+1/2\).
We set \(z=T^{-1/a}\). Then
\eqref{eq:varphi-def} and \eqref{eq:I-large} give
\begin{equation}\label{eq:analytic-relation}
 \frac rx=z^{-1}A(z^2),\qquad
 A(w)=(16/9)^{1/3}H_0
      \left(a+\sum_{m\ge1}c_mw^m\right).
\end{equation}
The function \(A\) is analytic near zero and satisfies \(A(0)>0\).
Writing \(\rho=x/r\), equation \eqref{eq:analytic-relation} becomes
\(\rho=z/A(z^2)\). The right-hand side is odd in \(z\) and has derivative
\(A(0)^{-1}>0\) at zero. The analytic inverse function theorem
therefore yields
\[
 z=\rho G(\rho^2),
\]
where \(G\) is analytic near zero and \(G(0)=A(0)>0\). Therefore
\begin{equation}\label{eq:F-analytic-extension}
 F=H_0x^az^{-a}
   =H_0r^aG\left((x/r)^2\right)^{-a}.
\end{equation}
Because \(G(0)>0\), the expression in
\eqref{eq:F-analytic-extension} is even and real analytic in
\(\xi\) near zero, and extends \(F\). At
\(\xi=0\), we have \(F=C_0r^a\), where
\(C_0=H_0G(0)^{-a}>0\), and hence
\[
 F_r=aC_0r^{a-1}>0,\qquad
 F_{rr}=a(a-1)C_0r^{a-2}>0.
\]
The identity \(4F_r\det D^2F=1\) extends by continuity. Thus
\(\det D^2F>0\), and the positivity of \(F_{rr}\) gives
\(D^2F>0\) at \(\xi=0\).

\smallskip
Finally, \(F/r^a\) is a positive continuous function of
\(\xi/r\in(-2,2)\). At the boundary, \(x/r\to2\) and
\(\varphi(r/x-1/2)\to H_0\), so
\[
 \frac{F}{r^a}\longrightarrow 2^aH_0
 \qquad\text{as }\frac{\xi}{r}\longrightarrow\pm2.
\]
Thus \(F/r^a\) extends to a positive continuous function on
\([-2,2]\). Its positive minimum and finite maximum prove
\(F\asymp r^a\) uniformly on \(\{r>|\xi|/2\}\).
\end{proof}

\begingroup
\setlength{\intextsep}{6pt}
\begin{figure}[H]
\raggedright
\noindent Figure~\ref{fig:F-profile} shows the normalized profile \(F/r^a\).
\par\smallskip
\centering
\begin{minipage}[c]{0.46\textwidth}
\centering
\begin{tikzpicture}
\begin{axis}[
 width=0.96\linewidth,
 height=0.64\linewidth,
 font=\small,
 tick label style={font=\footnotesize},
 xmin=-2.05,
 xmax=2.05,
 ymin=0.35,
 ymax=1.95,
 xlabel={\(q=\xi/r\)},
 ylabel={\(F(\xi,r)/r^{5/3}\)},
 xtick={-2,-1,0,1,2},
 ytick={0.5,1,1.5},
 axis lines=box,
 tick align=outside,
 grid=major,
 major grid style={draw=gray!20},
]
\addplot[blue!70!black,very thick,no marks]
 coordinates {
(-2.00,1.823399589997)
(-1.98,1.793772281239)
(-1.96,1.764892931997)
(-1.94,1.736568889760)
(-1.92,1.708731114821)
(-1.90,1.681339617131)
(-1.88,1.654367382040)
(-1.86,1.627794569733)
(-1.84,1.601605813238)
(-1.82,1.575788761250)
(-1.80,1.550333212015)
(-1.78,1.525230561223)
(-1.76,1.500473430925)
(-1.74,1.476055409520)
(-1.72,1.451970863346)
(-1.70,1.428214796310)
(-1.68,1.404782742848)
(-1.66,1.381670684657)
(-1.64,1.358874984809)
(-1.62,1.336392334827)
(-1.60,1.314219711617)
(-1.58,1.292354341997)
(-1.56,1.270793673195)
(-1.54,1.249535348066)
(-1.52,1.228577184102)
(-1.50,1.207917155519)
(-1.48,1.187553377861)
(-1.46,1.167484094677)
(-1.44,1.147707665930)
(-1.42,1.128222557858)
(-1.40,1.109027334049)
(-1.38,1.090120647563)
(-1.36,1.071501233936)
(-1.34,1.053167904944)
(-1.32,1.035119543024)
(-1.30,1.017355096268)
(-1.28,0.999873573904)
(-1.26,0.982674042217)
(-1.24,0.965755620845)
(-1.22,0.949117479413)
(-1.20,0.932758834459)
(-1.18,0.916678946620)
(-1.16,0.900877118055)
(-1.14,0.885352690074)
(-1.12,0.870105040947)
(-1.10,0.855133583882)
(-1.08,0.840437765154)
(-1.06,0.826017062360)
(-1.04,0.811870982803)
(-1.02,0.797999061978)
(-1.00,0.784400862159)
(-0.98,0.771075971072)
(-0.96,0.758024000653)
(-0.94,0.745244585873)
(-0.92,0.732737383635)
(-0.90,0.720502071732)
(-0.88,0.708538347848)
(-0.86,0.696845928627)
(-0.84,0.685424548777)
(-0.82,0.674273960213)
(-0.80,0.663393931249)
(-0.78,0.652784245818)
(-0.76,0.642444702723)
(-0.74,0.632375114922)
(-0.72,0.622575308838)
(-0.70,0.613045123696)
(-0.68,0.603784410879)
(-0.66,0.594793033311)
(-0.64,0.586070864854)
(-0.62,0.577617789732)
(-0.60,0.569433701963)
(-0.58,0.561518504818)
(-0.56,0.553872110284)
(-0.54,0.546494438553)
(-0.52,0.539385417517)
(-0.50,0.532544982287)
(-0.48,0.525973074711)
(-0.46,0.519669642922)
(-0.44,0.513634640887)
(-0.42,0.507868027975)
(-0.40,0.502369768538)
(-0.38,0.497139831504)
(-0.36,0.492178189985)
(-0.34,0.487484820898)
(-0.32,0.483059704606)
(-0.30,0.478902824564)
(-0.28,0.475014166988)
(-0.26,0.471393720541)
(-0.24,0.468041476032)
(-0.22,0.464957426660)
(-0.20,0.462141565083)
(-0.18,0.459593888528)
(-0.16,0.457314393194)
(-0.14,0.455303076731)
(-0.12,0.453559937518)
(-0.10,0.452084974492)
(-0.08,0.450878187002)
(-0.06,0.449939574691)
(-0.04,0.449269137377)
(-0.02,0.448866875001)
(0.00,0.448732787543)
(0.02,0.448866875001)
(0.04,0.449269137377)
(0.06,0.449939574691)
(0.08,0.450878187002)
(0.10,0.452084974492)
(0.12,0.453559937518)
(0.14,0.455303076731)
(0.16,0.457314393194)
(0.18,0.459593888528)
(0.20,0.462141565083)
(0.22,0.464957426660)
(0.24,0.468041476032)
(0.26,0.471393720541)
(0.28,0.475014166988)
(0.30,0.478902824564)
(0.32,0.483059704606)
(0.34,0.487484820898)
(0.36,0.492178189985)
(0.38,0.497139831504)
(0.40,0.502369768538)
(0.42,0.507868027975)
(0.44,0.513634640887)
(0.46,0.519669642922)
(0.48,0.525973074711)
(0.50,0.532544982287)
(0.52,0.539385417517)
(0.54,0.546494438553)
(0.56,0.553872110284)
(0.58,0.561518504818)
(0.60,0.569433701963)
(0.62,0.577617789732)
(0.64,0.586070864854)
(0.66,0.594793033311)
(0.68,0.603784410879)
(0.70,0.613045123696)
(0.72,0.622575308838)
(0.74,0.632375114922)
(0.76,0.642444702723)
(0.78,0.652784245818)
(0.80,0.663393931249)
(0.82,0.674273960213)
(0.84,0.685424548777)
(0.86,0.696845928627)
(0.88,0.708538347848)
(0.90,0.720502071732)
(0.92,0.732737383635)
(0.94,0.745244585873)
(0.96,0.758024000653)
(0.98,0.771075971072)
(1.00,0.784400862159)
(1.02,0.797999061978)
(1.04,0.811870982803)
(1.06,0.826017062360)
(1.08,0.840437765154)
(1.10,0.855133583882)
(1.12,0.870105040947)
(1.14,0.885352690074)
(1.16,0.900877118055)
(1.18,0.916678946620)
(1.20,0.932758834459)
(1.22,0.949117479413)
(1.24,0.965755620845)
(1.26,0.982674042217)
(1.28,0.999873573904)
(1.30,1.017355096268)
(1.32,1.035119543024)
(1.34,1.053167904944)
(1.36,1.071501233936)
(1.38,1.090120647563)
(1.40,1.109027334049)
(1.42,1.128222557858)
(1.44,1.147707665930)
(1.46,1.167484094677)
(1.48,1.187553377861)
(1.50,1.207917155519)
(1.52,1.228577184102)
(1.54,1.249535348066)
(1.56,1.270793673195)
(1.58,1.292354341997)
(1.60,1.314219711617)
(1.62,1.336392334827)
(1.64,1.358874984809)
(1.66,1.381670684657)
(1.68,1.404782742848)
(1.70,1.428214796310)
(1.72,1.451970863346)
(1.74,1.476055409520)
(1.76,1.500473430925)
(1.78,1.525230561223)
(1.80,1.550333212015)
(1.82,1.575788761250)
(1.84,1.601605813238)
(1.86,1.627794569733)
(1.88,1.654367382040)
(1.90,1.681339617131)
(1.92,1.708731114821)
(1.94,1.736568889760)
(1.96,1.764892931997)
(1.98,1.793772281239)
(2.00,1.823399589997)
};
\addplot[red!70!black,only marks,mark=o,
 mark options={fill=white,draw=red!70!black}]
 coordinates {(-2,1.823399589997) (2,1.823399589997)};
\end{axis}
\end{tikzpicture}
\end{minipage}\hfill
\begin{minipage}[c]{0.50\textwidth}
\setlength{\captionindent}{0pt}
\setlength{\abovecaptionskip}{0pt}
\raggedright
\caption[Normalized profile of \(F\).]{\raggedright
Normalized profile \(F(qr,r)/r^{5/3}\), \(-2<q<2\).
It is even and independent of \(r>0\). Open circles mark the
boundary value \(2^{5/3}H_0\).}
\label{fig:F-profile}
\end{minipage}
\end{figure}
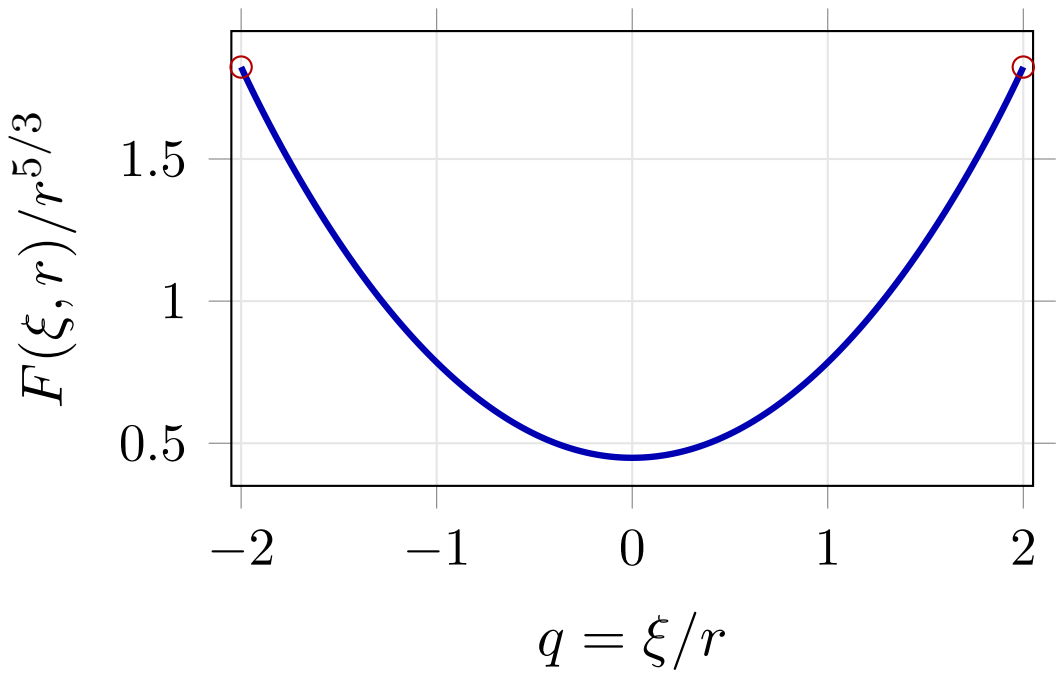
\endgroup

\medskip
We derive the boundary expansion of the outer potential to identify
the asymptotic terms that the smooth inner potential must reproduce
and to control the error when the two potentials are glued. We set
\[
 A_0=H_0a(a-1),\qquad B_0=(4A_0)^{-1/2},\qquad
 C_\Phi=\frac{2^{3/2}}3,\qquad
 D_0=\frac{(B_0C_\Phi)^2}{50H_0}.
\]
These are the constants used in \eqref{eq:model-outer-expansion};
in particular, \(A_0=10H_0/9\) and \(4A_0B_0^2=1\).
Writing \(h=T-1\), we have
\[
 \bigl((1+h)^{6/5}-1\bigr)^{-1/3}
 =(6h/5)^{-1/3}\bigl(1-h/30+O(h^2)\bigr).
\]
Integration in \eqref{eq:I-def} gives
\begin{equation}\label{eq:I-near-one}
 I(1+h)=\frac32(6/5)^{-1/3}h^{2/3}
       \left(1-\frac{h}{75}+O(h^2)\right).
\end{equation}
The factor in parentheses is analytic in \(h\). By
\eqref{eq:varphi-def}, raising this relation to the power \(3/2\)
and applying the analytic inverse function theorem gives
\begin{equation}\label{eq:T-near-one}
 T-1=\frac{y^{3/2}}{\sqrt5\,H_0^{3/2}}
       +\frac{y^3}{250H_0^3}+O(y^{9/2}),
\end{equation}
with a convergent power series in \(y^{3/2}\). Since
\[
 \frac1{\sqrt{5H_0}}=B_0C_\Phi,\qquad
 \frac1{250H_0^2}=D_0,
\]
substitution of \(y=s_0/x\), where \(s_0=r-x/2\), into
\eqref{eq:F} yields
\begin{equation}\label{eq:Fexp}
 F(\xi,r)=H_0x^a+B_0C_\Phi x^d s_0^p
       +D_0x^{d-p}s_0^{2p}
       +\wtO{x^{d-2p}s_0^{3p}}.
\end{equation}
This holds for \(0<s_0\le\eps x\), with fixed sufficiently small
\(\eps>0\). Here \(a-p=d\) and \(a-2p=d-p\); the weighted
derivative bounds follow by differentiating the convergent series
in \((s_0/x)^p\).

\subsection{The smooth fibre potential and its volume correction}

We construct the smooth fibre term \(B_0x^d\Phi(s)\) and the
correction \(B_0x^{d-p}\Psi(s)\) by computing their volume
density on \(Q(y)=1\).
We fix either base end, with \(x=|\log|u|^2|\gg1\), and choose
a local holomorphic square root of \(P(u)\). We set
\[
 \begin{aligned}
 z&=\sqrt P\,y,\qquad \sum_{j=1}^3y_j^2=1,\\
 \tau&=|y|^2=\frac{t}{|P(u)|}=\cosh s,\qquad
 s=s(u,z)=\arcosh\tau\ge0.
 \end{aligned}
\]
This identifies each fibre with the normalized quadric \(Q(y)=1\).
Let \(\zeta=\log u\) be a local base coordinate. Up to the fixed sign
determined by the residue convention, the residue form becomes
\begin{equation}\label{eq:residueproduct}
 \Om=\sqrt P\,\dd\log u\wedge\Om_Q .
\end{equation}
Here \(\Om_Q\) is the holomorphic residue two-form on \(Q(y)=1\).
The variables \(\tau\) and \(s\) are independent of the choice
of square root of \(P\).

\smallskip
For an invariant potential \(h=h(x,s)\), we set
\begin{equation}\label{eq:W}
 W(s)=\frac{\sinh s}{\sqrt{\cosh^2s+1}},
\end{equation}
and let
\(D^2h\) denote its real Hessian in \((x,s)\). By
\(\SO(3)\)-symmetry, for \(s>0\) we may compute at
\[
 y=(\cosh(s/2),\ii\sinh(s/2),0),
\]
using \((y_2,y_3)\) as local coordinates on the quadric. At this
point,
\[
 s_{y_2}=-\frac{\ii}{\cosh(s/2)},\qquad
 s_{y_2\bar y_2}=0,\qquad
 s_{y_3\bar y_3}=\frac1{\sinh s}.
\]
Since \(x=\pm(\zeta+\bar\zeta)\), the Hermitian determinant in
\((\zeta,y_2,y_3)\) is
\[
 \det(h_{j\bar k})
 =\frac{h_s\det D^2h}{\sinh s\cosh^2(s/2)}.
\]
Using \(\Om_Q=\dd y_2\wedge\dd y_3/(2y_1)\) and
\eqref{eq:residueproduct}, we obtain
\[
 \frac{(\ddc h)^3}{6(\ii)^9\Om\wedge\bar\Om}
 =\frac{4h_s\det D^2h}{|P|\sinh s}.
\]
Since \(t=|P|\cosh s\), dividing by the model volume
form \(\nu_i\) from Subsection~\ref{subsec:potentials} gives
\begin{equation}\label{eq:invariant-volume}
 \mathcal V(h):=\frac{(\ddc h)^3}{\nu_i}
 =\frac{4h_s\det D^2h}{W}.
\end{equation}
If \(h\) is smooth in \((x,\tau)\), then \(\mathcal V(h)\) extends
smoothly across \(s=0\).

\medskip
Primes denote derivatives with respect to \(s\). For the leading
potential \(h_0=H_0x^a+B_0x^d\Phi(s)\), we write \(q=\Phi'\).
Its leading normalized volume density is
\[
 4A_0B_0^2x^{a+2d-2}\frac{qq'}W=\frac{qq'}W,
\]
since \(a+2d=2\) and \(4A_0B_0^2=1\).
Thus the leading volume equation is \(qq'=W\), with
\(q(0)=0\) for smoothness at \(s=0\). Since
\(W=\partial_s\arsinh(\cosh s)\), integration gives
\begin{equation}\label{eq:fibre}
 \Phi'=q=\sqrt{2\{\arsinh(\cosh s)-\arsinh1\}}.
\end{equation}
Near \(s=0\), \(q/s\) is a smooth positive function of \(s^2\).
Consequently, \(\Phi\) is smooth as a function of
\(\tau=\cosh s\), including at the vanishing sphere \(\tau=1\).
The fibre Hessian coefficients \(q'\) and \(q/\sinh s\) are positive
for \(s>0\) and extend to positive smooth functions at \(s=0\). Hence
\(\ddc\Phi\) defines a K\"ahler metric \(g_Q\) on the entire
quadric. Moreover,
\[
 q\sim\sqrt{2s},\qquad q'\sim(2s)^{-1/2}
 \qquad (s\to\infty).
\]
The \(\SO(3)\)-orbits are compact, and the radial length element is
a fixed multiple of \(\sqrt{q'}\,\dd s\). Since
\(\int^\infty\sqrt{q'}\,\dd s=\infty\), the metric \(g_Q\) is complete.

\smallskip
For large \(s\), \(\arsinh(\cosh s)=s+O(e^{-2s})\), hence
\[
 q=\sqrt{2(s-\arsinh1)}
   +O\!\left(s^{-1/2}e^{-2s}\right).
\]
The remainder is integrable at infinity. We therefore choose the
additive constant in \(\Phi\) so that
\(\Phi-C_\Phi(s-\arsinh1)^p\to0\).
Using \(q'=W/q\) and \(W=1+O(e^{-2s})\), we obtain
\begin{equation}\label{eq:gexp}
 \Phi(s)=C_\Phi(s-\arsinh1)^p+O(s^2e^{-2s}),\qquad
 S(s):=\frac{q^2}{q'}=3\Phi(s)+O(s^2e^{-2s}).
\end{equation}
Each remainder and its \(j\)-th derivative are bounded by
\(C_js^2e^{-2s}\), for every fixed \(j\ge0\).
Expressions involving \((s-\arsinh1)^p\) are used only for
large \(s\).

\medskip
For \(h_0=H_0x^a+B_0x^d\Phi\), direct differentiation gives
\[
 \begin{aligned}
 (h_0)_s&=B_0x^dq,\\
 \det D^2h_0
 &=A_0B_0x^{a+d-2}q'
   +B_0^2x^{2d-2}
      \{d(d-1)\Phi q'-d^2q^2\}.
 \end{aligned}
\]
Using \(4A_0B_0^2=1\), \(a+2d=2\), and \(3d-2=-p\), we obtain
\begin{equation}\label{eq:uncorrected-error}
 \mathcal V(h_0)=1+x^{-p}E_1,\qquad
 E_1=\frac{B_0}{A_0}\{d(d-1)\Phi-d^2S\}.
\end{equation}
The radial linearized operator is
\[
 \mathcal L_Qv=\frac{v'}q+\frac{v''}{q'}
              =\frac{(qv')'}W.
\]
Adding \(B_0x^{d-p}\Psi\) changes the coefficient of \(x^{-p}\) in
the normalized density to \(E_1+\mathcal L_Q\Psi\). To cancel it,
we impose \(\mathcal L_Q\Psi=-E_1\). Since \(d=1/6\) and
\(A_0=10H_0/9\), this gives
\begin{equation}\label{eq:model-correction-equation}
 \frac{\Psi'}q+\frac{\Psi''}{q'}
 =\frac{B_0}{40H_0}\left(5\Phi+\frac{q^2}{q'}\right),
 \qquad s>0.
\end{equation}
The regular solution satisfies
\begin{equation}\label{eq:corrector}
 (q\Psi')'=-WE_1,\qquad
 \lim_{s\downarrow0}q(s)\Psi'(s)=0.
\end{equation}
The boundary condition gives
\[
 \Psi'(s)=-\frac1{q(s)}\int_0^s W(v)E_1(v)\,\dd v.
\]
The integral formula determines \(\Psi\) up to an additive constant.
Near zero,
\(WE_1\) is \(s\) times a smooth function of \(s^2\), while
\(q/s\) is smooth and positive. The integral formula shows that
\(\Psi'/s\) is smooth in \(s^2\), so \(\Psi\) is smooth in
\(\tau=\cosh s\) at \(\tau=1\).

\begin{lemma}\label{lem:corrector}
We set \(c_2=\frac{B_0}{A_0}\frac{d(1+2d)}{10}
          =\frac{B_0}{50H_0}\).
There is a unique choice of the additive constant in \(\Psi\) for
which, as \(s\to\infty\),
\begin{equation}\label{eq:psiexp}
 \Psi(s)=c_2\Phi(s)^2+c_\infty q(s)+O(s^2e^{-2s})
\end{equation}
for a constant \(c_\infty\). The remainder satisfies the derivative
bounds in \eqref{eq:gexp}, and \(B_0c_2C_\Phi^2=D_0\).
\end{lemma}
\begin{proof}
Using \(\Phi'=q\) and \(S=q^2/q'\), we compute
\[
 \mathcal L_Q(\Phi^2)=4\Phi+2S,\qquad
 E_1+c_2\mathcal L_Q(\Phi^2)
 =\frac{B_0}{A_0}\frac{d(1-3d)}5(S-3\Phi).
\]
Thus \(U=\Psi-c_2\Phi^2\) satisfies, by \eqref{eq:corrector},
\[
 (qU')'=-W\bigl(E_1+c_2\mathcal L_Q(\Phi^2)\bigr).
\]
By \eqref{eq:gexp}, the right-hand side is exponentially small,
and integration gives
\[
 qU'=c_\infty+O(s^2e^{-2s})
\]
for a finite constant \(c_\infty\). Since
\(q'=W/q\) and \(W=1+O(e^{-2s})\), a second integration gives
\[
 U=c_\infty q+c_0+O(s^2e^{-2s}),
\]
with the same derivative bounds. Subtracting \(c_0\) gives
\eqref{eq:psiexp}. This normalization is unique: a change in
\(c_\infty\) contributes an unbounded multiple of \(q\), and a
nonzero additive constant cannot be absorbed by the remainder.

The identities \(5d(1+2d)=a(a-1)\) and
\(A_0=H_0a(a-1)\), with the definition of \(c_2\), give
\[
 B_0c_2C_\Phi^2
 =\frac{B_0^2C_\Phi^2}{50H_0}=D_0.
\]
\end{proof}

\medskip
The inner potential is therefore
\begin{equation}\label{eq:Ki}
 K_i(u,z)=H_0x^a+B_0x^d\Phi(s)+B_0x^{d-p}\Psi(s).
\end{equation}
It is independent of the local choice of \(\sqrt P\).
By \(t=|P|\cosh s\) and \(\log|P|=x/2+O(e^{-x})\), the
outer variable \(s_0=r-x/2\) satisfies
\[
 s_0=s-\arsinh1+O(e^{-x})+O(e^{-2s})
 \qquad (x,s\to\infty).
\]
By \eqref{eq:gexp}, \eqref{eq:psiexp}, and
\(B_0c_2C_\Phi^2=D_0\), the leading fibre term and its quadratic
correction agree with those in \eqref{eq:Fexp}. The remaining term
\(B_0c_\infty x^{d-p}q\) has order \(x^{d-p}s^{1/2}\); it is
included in the gluing estimate \eqref{eq:deltader}.

\begin{lemma}\label{lem:inner}
For sufficiently large \(x\), the following hold on
\(0\le s\le3x^{1/4}\).
\begin{enumerate}
\item[(i)] The potential \(K_i\) is smooth and strictly
plurisubharmonic. Its associated Riemannian metric is uniformly
comparable to
\[
 g_*=2A_0x^{a-2}|\dd\zeta|^2+B_0x^d g_Q,
\]
where \(g_Q\) is the metric associated with \(\ddc\Phi\), pulled
back from \(Q(y)=1\) by \(y=z/\sqrt P\).

\item[(ii)] The normalized volume density satisfies
\begin{equation}\label{eq:innererror}
 \frac{4(K_i)_s\det D^2K_i}{W}
 =1+\wtO{\frac{(1+s)^3}{x^3}}.
\end{equation}
The normalization is that of \eqref{eq:invariant-volume}. At \(s=0\),
the quotient denotes its smooth extension in
\(\tau=\cosh s\).
\end{enumerate}
\end{lemma}
\begin{proof}
Since \(\Phi\) and \(\Psi\) are smooth in \(\tau\), so is \(K_i\).
We set \(\eta=x^{-p}(1+s)^p\). Differentiating \eqref{eq:Ki} and using
\eqref{eq:gexp} and \eqref{eq:psiexp} gives
\[
 \begin{aligned}
 (K_i)_{xx}&=A_0x^{a-2}(1+O(\eta)),\\
 (K_i)_{ss}&=B_0x^d q'(1+O(\eta)),\qquad
 (K_i)_s=B_0x^d q(1+O(\eta)),\\
 \frac{|(K_i)_{xs}|}{\sqrt{A_0x^{a-2}\,B_0x^d q'}}
 &\le C\eta^{1/2}.
 \end{aligned}
\]
The estimates hold near \(s=0\) by smoothness in \(\tau\); in
particular, \(\Psi'/q\) and \(\Psi''/q'\) extend smoothly there.
Since \(q'>0\) and \(q/\sinh s>0\), the diagonal blocks of the Levi
form are positive and comparable to those of \(g_*\).
On the stated region, \(\eta=O(x^{-9/8})\), so the normalized mixed
block tends uniformly to zero. The Schur complement is therefore
positive for large \(x\), proving (i), including at \(s=0\).

\smallskip
For the volume estimate, we set \(\lambda=x^{-p}\).
With \(a-d=p\) and \(qq'=W\), the determinant in
\eqref{eq:invariant-volume} factors as
\[
\begin{aligned}
 \mathcal V(K_i)
 &=\left(1+\lambda\frac{\Psi'}q\right)
 \Biggl[\left(1+\lambda\frac{\Psi''}{q'}\right)
 \left(1+\frac{B_0}{A_0}d(d-1)\lambda\Phi
 +\frac{B_0}{A_0}(d-p)(d-p-1)\lambda^2\Psi\right)\\
 &\hspace{28mm}
 -\frac{B_0}{A_0}\lambda\frac{q^2}{q'}
   \left(d+(d-p)\lambda\frac{\Psi'}q\right)^2\Biggr].
\end{aligned}
\]
The coefficient of \(\lambda\) is
\[
 \frac{\Psi'}q+\frac{\Psi''}{q'}
 +\frac{B_0}{A_0}\left(d(d-1)\Phi-d^2\frac{q^2}{q'}\right)
 =\mathcal L_Q\Psi+E_1=0
\]
by \eqref{eq:corrector}. The remaining terms have degrees two
through four in \(\lambda\). By \eqref{eq:gexp},
\eqref{eq:psiexp}, and smoothness at zero,
\[
 \Phi=O((1+s)^p),\quad
 \Psi=O((1+s)^{2p}),\quad
 \frac{q^2}{q'},\ \frac{\Psi'}q,\ \frac{\Psi''}{q'}
       =O((1+s)^p),
\]
with weighted derivatives. Since
\(\lambda(1+s)^p=O(x^{-9/8})\), the terms of degree \(j\ge2\)
are \(\wtO{\lambda^j(1+s)^{jp}}\), and their sum is
\[
 \mathcal V(K_i)-1
 =\wtO{\lambda^2(1+s)^{2p}}
 =\wtO{\frac{(1+s)^3}{x^3}}.\qedhere
\]
\end{proof}

\subsection{Estimates for the glued potential}\label{app:gluing-volume}

For \(\SO(3)\)-invariant functions, \(O_{\mathrm{wt}}\) denotes
bounds for all fixed-order weighted derivatives. We write
\(\theta=\arg u\). On either base end, for \(s\ge1\),
\(f=\wtO{x^m s^n}\) means
\[
 |\partial_x^i\partial_s^j\partial_\theta^k f|
 \le C_{ijk}x^{m-i}s^{n-j},\qquad i,j,k\in\Z_{\ge0}.
\]
Near \(s=0\), we replace \(s\) in the weights by \(1+s\) and
take fibre derivatives in smooth local coordinates on \(Q(y)=1\).

\smallskip
On a central subcone \(|\xi|\le(2-\delta)\log t\), where
\(\xi=\log|u|^2\) and \(0<\delta<2\) is fixed, we use
\((\log t)\partial_\xi\), \((\log t)t\partial_t\), and
\(\partial_\theta\). Equivalently,
\(f=\wtO{(\log t)^m}\) means
\[
 |\partial_\xi^i(t\partial_t)^j\partial_\theta^k f|
 \le C_{ijk,\delta}(\log t)^{m-i-j},\qquad i,j,k\in\Z_{\ge0}.
\]
Constants may depend on the derivative orders, the fixed coordinate
domains, and \(\delta\), but are uniform in the angular variables.

\begin{proof}[Proof of \eqref{eq:gluing-metric-error}]
Invariance under
\(\SO(3)\) and \((u,z)\mapsto(u^{-1},z)\) reduces the calculation
to \(|u|>1\), where
\[
 \log|P|-x/2=O(e^{-x}),\qquad
 s-\log(2t/|P|)=O(e^{-2s}),
\]
with the same bounds after any fixed number of derivatives,
uniformly in \(\arg u\). By \eqref{eq:Fexp},
\eqref{eq:gexp}, and \eqref{eq:psiexp}, the common base term
and the first two fibre terms cancel in \(K_o-K_i\). Hence, on
the transition strip,
\begin{equation}\label{eq:deltader}
 |\partial_x^i\partial_s^j(K_o-K_i)|
 \le\frac{C_{ij}}{x^i s^j}
 \left(\frac{s^{9/2}}{x^{17/6}}+\frac{s^{1/2}}{x^{4/3}}\right)
 +\frac{C_{ijN}}{x^N}
\end{equation}
for all fixed nonnegative integers \(i,j\) and \(N>0\).
The first term in parentheses comes from the remainder in
\eqref{eq:Fexp}; the second comes from \(-B_0c_\infty x^{-4/3}\Phi'(s)\). The other terms and their derivatives, including
angular derivatives, are exponentially small in \(x\) or \(s\).
Since \(s\asymp x^{1/4}\), they are bounded by \(C_N/x^N\)
for every \(N>0\).

\smallskip
The derivative bounds in \eqref{eq:gexp} give
\(\Phi'\asymp s^{1/2}\) and \(\Phi''\asymp s^{-1/2}\) for
large \(s\). With \(g_*\) as in \eqref{eq:model-metric}, an
invariant function \(h=h(x,s)\) therefore satisfies
\[
 \|\ddc h\|_{g_*}\le C\Bigl(x^{1/3}|h_{xx}|
       +x^{1/12}s^{1/4}|h_{xs}|
       +\frac{s^{1/2}}{x^{1/6}}|h_{ss}|
       +\frac{1}{x^{1/6}s^{1/2}}|h_s|\Bigr).
\]
These are the base, mixed, radial, and transverse contributions
in the local coordinates \((\log u,y)\), with
\(y=z/\sqrt P\). Applying this bound to \eqref{eq:deltader}
and absorbing the angular remainder gives
\begin{equation}\label{eq:deltalevi}
 \|\ddc(K_o-K_i)\|_{g_*}
 \le C\left(\frac{s^3}{x^3}+\frac1{x^{3/2}s}\right)
 +\frac{C_N}{x^N}.
\end{equation}

\smallskip
The cutoff satisfies
\begin{equation}\label{eq:cutoff-bounds}
 \bigl|\partial_x^i\partial_s^j\vartheta(s/x^{1/4})\bigr|
 \le\frac{C_{ij}}{x^i s^j}
 \qquad\text{on }x^{1/4}\le s\le2x^{1/4}.
\end{equation}
The product rule therefore preserves the orders in
\eqref{eq:deltader} for
\(K_\infty-K_i=\vartheta(s/x^{1/4})(K_o-K_i)\).
Applying \eqref{eq:deltalevi} with \(N\ge2\) and
\(s\asymp x^{1/4}\) gives \eqref{eq:gluing-metric-error}.
\end{proof}

\begin{proof}[Proof of Lemma~\ref{lem:volumeerror}]
Outside a compact set, \(\omega_0=\ddc K_\infty\).
On either end \(x>x_0\), the identity
\(|u-u^{-1}|^2=|P(u)|^2-4\cos(2\theta)\) gives
\[
 \frac{\nu_i}{\nu}
 =\left(1+\frac{2t+1-4\cos(2\theta)}{t^2+|P|^2}\right)^{1/2}
 =1+O\!\left(\frac1t\right).
\]
For \(0\le s\le2x^{1/4}\), we have
\(t\ge|P|\asymp e^{x/2}\), so \(\nu_i/\nu-1\) and its
weighted derivatives are exponentially small in \(x\).
Multiplying by the estimate in Lemma~\ref{lem:inner} yields
\begin{equation}\label{eq:inner-volume-error}
 \frac{(\ddc K_i)^3}{\nu}
 =1+\wtO{\frac{(1+s)^3}{x^3}}
 =1+\wtO{\frac1{x^{9/4}}}.
\end{equation}
On \(0\le s\le x^{1/4}\), we have
\(K_\infty=K_i\), so \eqref{eq:inner-volume-error} also
estimates \(\omega_0^3/\nu\).

\smallskip
On \(x^{1/4}\le s\le2x^{1/4}\),
\eqref{eq:gluing-metric-error} bounds the metric perturbation by
\(C/x^{7/4}\). Equations~\eqref{eq:deltader} and
\eqref{eq:cutoff-bounds} give the same order for all fixed-order
weighted derivatives of its normalized coefficients.
By Lemma~\ref{lem:inner},
the inner metric is uniformly comparable to \(g_*\), and the
normalized metric matrix and its inverse have bounded weighted
derivatives. Thus
\begin{equation}\label{eq:transition-volume-ratio}
 \frac{(\ddc K_\infty)^3}{(\ddc K_i)^3}
 =1+\wtO{\frac1{x^{7/4}}}.
\end{equation}
Combining \eqref{eq:inner-volume-error} and
\eqref{eq:transition-volume-ratio}, we obtain
\[
 \frac{\omega_0^3}{\nu}
 =\frac{(\ddc K_\infty)^3}{(\ddc K_i)^3}
   \frac{(\ddc K_i)^3}{\nu}
 =1+\wtO{\frac1{x^{7/4}}}.
\]

\smallskip
Where \(K_\infty=K_o\), we use ambient coordinates
\((\zeta,z)\), with \(\zeta=\log u\). Let \(g_{\rm amb}\)
be the complex Hessian of \(K_o\). We have \(\partial\xi=\dd\zeta\)
and, since \(r=\log t+\log 2-\arsinh1\),
\[
 \partial r=\frac{1}{t}\sum_j\bar z_j\,\dd z_j,\qquad
 r_{z_j\bar z_k}=\frac{\delta_{jk}}{t}-\frac{\bar z_jz_k}{t^2}.
\]
In a unitary frame of the \(z\)-space with first vector \(z/\sqrt t\),
\(\partial r\) has radial component \(t^{-1/2}\) and vanishes on the
other two vectors. The complex Hessian of \(r\) vanishes in the radial
direction and is \(t^{-1}I_2\) on its orthogonal complement. The chain
rule gives
\[
 g_{\rm amb}
 =\begin{pmatrix}
 F_{\xi\xi} & F_{\xi r}/\sqrt t\\
 F_{\xi r}/\sqrt t & F_{rr}/t
 \end{pmatrix}
 \oplus\frac{F_r}{t}I_2.
\]
The first block is obtained from \(D^2F\) by congruence with
\(\operatorname{diag}(1,t^{-1/2})\). It has determinant
\(t^{-1}\det D^2F\), so \(\det g_{\rm amb}=F_r^2\det D^2F/t^3\).

For the defining function \(C=Q-P\), the base component of
\(\partial C\) is \(-(u-u^{-1})\). On \(X\), the radial component is
\(\dd Q(z/\sqrt t)=2P/\sqrt t\), and
\(\lvert\dd Q\rvert_{\rm Eucl}^2=4t\). The transverse Euclidean
squared norm is therefore \(4t-4|P|^2/t\), which is multiplied by
\(t/F_r\) in the inverse metric. For the base and radial contribution,
inverting the first block multiplies the radial covector component
by \(\sqrt t\) before taking the squared norm with \((D^2F)^{-1}\).
Thus the squared conormal norm is
\[
 N_C=|\partial C|_{g_{\rm amb}}^2
 =\frac{4(t^2-|P|^2)}{F_r}
 + (-(u-u^{-1}),2P)(D^2F)^{-1}
 (-\overline{u-u^{-1}},2\bar P)^T.
\]
To estimate the inverse Hessian, we set \(M=F_r(D^2F)^{-1}\).
Equation~\eqref{eq:model-outer-equation} gives
\[
 \det D^2F=\frac1{4F_r},\qquad
 M=4F_r^2\operatorname{adj}(D^2F).
\]
Near either cone boundary, we set \(x=|\xi|\) and \(s_0=r-x/2\).
For \(0<s_0\le\eps x\), the differentiated expansion
\eqref{eq:Fexp} gives, in the \((\xi,r)\)-coordinates,
\[
 F_r=O(x^{1/6}s_0^{1/2}),\qquad
 \|D^2F\|\le C\bigl(x^{-1/3}+x^{1/6}s_0^{-1/2}\bigr).
\]
Since the adjugate of a two-by-two matrix has comparable norm,
\[
 \|M\|\le C\bigl(s_0+\sqrt{xs_0}\bigr)\le Cr.
\]
On the remaining fixed interior subcone, \(M\) is homogeneous
of degree one and its profile on \(r=1\) is smooth and bounded.
Thus the same bound holds there, and throughout the outer domain
we obtain \(\|F_r(D^2F)^{-1}\|\le Cr\le C\log t\).
The adjugate identity, with the residue normalization of \(\Om\),
gives the first equality below; \eqref{eq:model-outer-equation}
and \(|u-u^{-1}|^2\le|P|^2+4\) give the second:
\[
 \frac{(\ddc K_o)^3}{6(\ii)^9\Om\wedge\bar\Om}
 =(\det g_{\rm amb})N_C
 =\frac1t+O\!\left(\frac{\log t\,(|P|^2+1)}{t^3}\right).
\]
Dividing by \(\nu\) from \eqref{eq:nu} gives
\begin{equation}\label{eq:outer-volume-error}
 \frac{(\ddc K_o)^3}{\nu}
 =1+O\!\left(\frac1t+\frac{\log t\,(|P|^2+1)}{t^2}\right).
\end{equation}

\smallskip
On either end, \(K_\infty=K_o\) for \(s\ge2x^{1/4}\), where
\[
 \frac{|P|}{t}=\frac1{\cosh s},\qquad
 \log t=\frac x2+\log\cosh s+O(e^{-x}).
\]
Since \(x\le(s/2)^4\), we have \(\log t\le C(1+s^4)\) and
\[
 \frac1t+\frac{\log t\,(|P|^2+1)}{t^2}
 \le C\bigl(e^{-s}+(1+s^4)e^{-2s}\bigr).
\]
Weighted differentiation introduces only polynomial factors
in \(x,s\). Thus \eqref{eq:outer-volume-error} has remainder
\(\wtO{1/(\log t)^N}\) for every \(N>0\).

\smallskip
On a central subcone \(|\xi|\le(2-\delta)\log t\),
\(0<\delta<2\), the bound \(|P(u)|\le2e^{|\xi|/2}\) gives
\[
 \frac{|P|^2}{t^2}\le\frac4{t^\delta},\qquad
 \frac1t+\frac{\log t\,(|P|^2+1)}{t^2}
 \le C_\delta\left(\frac1t+\frac{\log t}{t^\delta}\right).
\]
For large \(t\), we have \(K_\infty=K_o\) throughout this
subcone: on the ends, the first inequality gives
\(s\ge(\delta/2)\log t-\log2>2x^{1/4}\), while the central
prescription covers bounded \(x\), including \(P(u)=0\).
The weighted derivatives introduce only powers of
\(\log t\), so \eqref{eq:outer-volume-error} again has
remainder \(\wtO{1/(\log t)^N}\) for every \(N>0\).
All estimates include \(\theta\)-derivatives; the other angular
variables are absent by \(\SO(3)\)-invariance.

\smallskip
Since \(\log t=x/2+O(x^{1/4})\asymp x\) on the inner and
transition regions, the inner, transition, and outer estimates prove
\eqref{eq:volumeerror}.
\end{proof}

\subsection{The Hessian metric as a flat sector}

\begin{lemma}\label{lem:hessiancone}
The metric \(g_H=\tfrac12D^2F\) on \(\{r>|\xi|/2\}\) is
isometric to the open Euclidean sector of opening angle \(5\pi/6\). Its
radial coordinate is
\[
 R_F=\sqrt{\frac{4F}{5}}.
\]
On the half-cone \(\xi>0\), the angular coordinate measured from
the boundary ray \(r=\xi/2\) is
\begin{equation}\label{eq:sigma-coordinate}
 \sigma(T)=\int_1^T\frac{\dd\tau}
                 {2\tau\sqrt{\tau^{6/5}-1}},
 \qquad T=\frac{\varphi(y)}{H_0},\qquad
 y=\frac r\xi-\frac12.
\end{equation}
Its range is \((0,5\pi/12)\).
\end{lemma}
\begin{proof}
We first work on the half-cone \(\xi=x>0\). The variables
\[
 x>0,\qquad y=\frac rx-\frac12>0
\]
are global coordinates there. By \eqref{eq:varphi-def},
\(T=\varphi(y)/H_0\) increases from \(1\) to \(\infty\). Since
\(F=x^aH_0T\), the pair \((F,T)\) also gives global
coordinates, with range \((0,\infty)\times(1,\infty)\).

\smallskip
Substituting \(\dd r=(y+1/2)\dd x+x\dd y\) into
\eqref{eq:F-second-derivatives} gives
\[
 2g_H=a(a-1)x^{a-2}\varphi\,(\dd x)^2
 +2(a-1)x^{a-1}\varphi'\,\dd x\,\dd y
 +x^a\varphi''\,(\dd y)^2.
\]
Since \(\dd F=ax^{a-1}\varphi\,\dd x+x^a\varphi'\,\dd y\),
completing the square yields
\begin{equation}\label{eq:hessian-polar-form}
 g_H=\frac{a-1}{2aF}(\dd F)^2
 +\frac F2\left(
       \frac{\varphi''}{\varphi}
       -\frac{a-1}{a}\frac{(\varphi')^2}{\varphi^2}
       \right)(\dd y)^2.
\end{equation}
Because \(a=5/3\), we have
\(R_F^2=2(a-1)F/a=4F/5\), so the first term in
\eqref{eq:hessian-polar-form} is \((\dd R_F)^2\). Dividing
the second term by \(R_F^2\) gives
\[
 \frac{a\varphi\varphi''-(a-1)(\varphi')^2}
      {4(a-1)\varphi^2}(\dd y)^2.
\]
Substituting \eqref{eq:varphi-identity},
\(\varphi=H_0T\), \(\dd y=H_0\dd T/\varphi'\), and the first
identity in \eqref{eq:varphi-ode} reduces this expression to
\[
 \frac{(\dd T)^2}{4T^2(T^{6/5}-1)}=(\dd\sigma)^2,
\]
where \(\sigma\) is defined by \eqref{eq:sigma-coordinate}.
Consequently,
\[
 g_H=(\dd R_F)^2+R_F^2(\dd\sigma)^2
\]
on \(\xi>0\), so this half-cone is flat. Its angular width is
\begin{equation}\label{eq:half-sector-angle}
 \int_1^\infty\frac{\dd T}{2T\sqrt{T^{6/5}-1}}
 =\frac56\int_0^1\frac{\dd u}{\sqrt{1-u^2}}
 =\frac{5\pi}{12},
\end{equation}
where we used \(u=T^{-3/5}\).

\smallskip
At fixed \(F\), the points approach the boundary ray \(r=\xi/2\)
as \(T\downarrow1\), and approach \(\xi=0\) as \(T\to\infty\).
Lemma~\ref{lem:F} shows that the even extension of \(F\) is smooth
across \(\xi=0\). The half-cone \(\xi<0\) is its reflected copy,
again of angular width \(5\pi/12\). Gluing the two copies along
\(\xi=0\) therefore gives a flat sector of total opening angle
\(5\pi/6\).
\end{proof}

\section{Holomorphic charts and the covariant Hessian bound}\label{app:atlas}

We use the model and gluing estimates in Section~\ref{sec:background}
and Appendix~\ref{app:model} to prove Proposition~\ref{prop:atlas}
and the bound \(|\nabla^2K_0|_{g_0}\le C\) needed in
Proposition~\ref{prop:complete}.

\begin{proof}[Proof of Proposition~\ref{prop:atlas}]
Let \(B\subset\C^3\) be the unit ball. We construct holomorphic
local diffeomorphisms \(\Phi_p:4B\to X\) with \(\Phi_p(0)=p\),
estimate the pullback metrics on \(3B\), and restrict to \(2B\).
All scales are evaluated at \(p\), and the constant \(\kappa>0\)
is fixed and independent of \(p\).

\smallskip
We choose \(S\gg1\) and then \(\eps_0>0\) sufficiently small that
\eqref{eq:Fexp} applies on \(S/2\le s\le3\eps_0x\).
Suppose first that \(S\le s(p)\le2\eps_0x(p)\), where
\(x=\xi\) on the positive end and \(x=-\xi\) on the negative end.
We choose a logarithm \(\zeta_p\) of \(u(p)\) and a holomorphic
square root of \(P(e^\zeta)\) on the \(\zeta\)-domain.
After an \(\SO(3)\) rotation, we use the coordinates
\begin{equation}\label{eq:endchart}
 \begin{gathered}
 u=e^\zeta,\qquad z=\sqrt P\,y,\qquad
 y_3=e^\lambda v,\\
 y_1-\ii y_2=e^\lambda,\qquad
 y_1+\ii y_2=e^{-\lambda}-e^\lambda v^2.
 \end{gathered}
\end{equation}
They satisfy \(Q(y)=1\), and the centre is
\((\zeta,\lambda,v)=(\zeta_p,s(p)/2,0)\).
We set
\begin{equation}\label{eq:steps}
 l_x=x^{1/6},\qquad l_s=x^{-1/12}s^{1/4},\qquad
 l_v=x^{-1/12}s^{-1/4},
\end{equation}
and define \(\Phi_p\) by substituting
\[
 \zeta=\zeta_p+\kappa l_x\eta_1,\qquad
 \lambda=s(p)/2+\kappa l_s\eta_2,\qquad
 v=\kappa l_v\eta_3,\qquad \eta\in4B,
\]
in \eqref{eq:endchart}. For large \(x(p)\), this domain stays
over the same base end, where \(P\ne0\).

\smallskip
Let \(h=h(x,s)\) be obtained from \(K_\infty\) by replacing
\(\log|P|\) with \(x/2\). Its outer term is
\(F(x,x/2+\log(2\cosh s)-\arsinh1)\).
With \(\theta=\arg u\),
\[
 \log|P|-\frac x2
 =\operatorname{Re}\log(1+e^{-x\mp2\ii\theta})=O(e^{-x}),
\]
with the same bound after fixed-order derivatives in \(x,\theta\).
Thus \(K_\infty-h\) and its scaled derivatives are bounded by
\(C_j e^{-x(p)/2}\). The expansions \eqref{eq:Fexp},
\eqref{eq:gexp}, and \eqref{eq:psiexp}, together with the cutoff
bounds in Lemma~\ref{lem:gluing}, give
\[
 \begin{aligned}
 |\partial_x^i\partial_s^jh|
 &\le C_{ij}x^{1/6-i}s^{3/2-j} &&(j\ge1),\\
 |\partial_x^i(h-H_0x^{5/3})|
 &\le C_i x^{1/6-i}s^{3/2}.&
 \end{aligned}
\]
At the centre, the diagonal Levi coefficients have orders
\[
 h_{xx}\asymp x^{-1/3},\qquad
 h_{ss}\asymp x^{1/6}s^{-1/2},\qquad
 h_s\frac{e^s}{\sinh s}\asymp x^{1/6}s^{1/2}.
\]
The scales in \eqref{eq:steps} therefore give
\[
 l_x^2h_{xx}\asymp l_s^2h_{ss}
 \asymp l_v^2h_s\frac{e^s}{\sinh s}\asymp1,
 \qquad |l_xl_sh_{xs}|\le C(s/x)^{3/4}.
\]
For sufficiently small \(\eps_0\), the normalized mixed entry
is smaller than the diagonal lower bound. The metric is then
uniformly positive at the centre.

\smallskip
For \(i+j\ge2\) and \(j\ge1\),
\[
 l_x^il_s^j|\partial_x^i\partial_s^jh|
 \le C_{ij}x^{1/6-5i/6-j/12}s^{3/2-3j/4}.
\]
Since \(S/2\le s\le3\eps_0x\) throughout the chart, the powers
on the right are bounded independently of \(p\):
\begin{equation}\label{eq:exponent}
 \frac16-\frac{5i}6-\frac j{12}
       +\max\left\{\frac32-\frac{3j}4,0\right\}\le0.
\end{equation}
For \(j=1\), the left side is \(5(1-i)/6\le0\); for \(j=2\),
it is \(-5i/6\); and for \(j\ge3\), it is at most
\(1/6-j/12<0\). Pure base derivatives satisfy
\(l_x^i|\partial_x^ih|\le C_i x^{5/3-5i/6}\le C_i\) for \(i\ge2\).

\smallskip
For the nonlinear terms in the chain rule, \eqref{eq:endchart}
gives the exact identity
\[
 \tau=\frac{e^{\lambda+\bar\lambda}}2
       \left(1+2|v|^2+|e^{-2\lambda}-v^2|^2\right),
\]
so
\[
 s=\lambda+\bar\lambda+2\log(1+|v|^2)+\mathcal R,
\]
where \(\mathcal R\) and its fixed-order scaled derivatives
are bounded by a polynomial in \(s\) times \(e^{-s}\).
Since \(v=O(l_v)\),
every positive-order derivative of \(2\log(1+|v|^2)\) in the
\(\eta\)-coordinates is \(O(l_v^2)\). The identities
\begin{equation}\label{eq:nonlinear}
 h_s l_v^2\asymp1,\qquad
 \frac{l_v^2}{l_s}=\frac{l_s}{s}
     =x^{-1/12}s^{-3/4}\le1,\qquad h_s l_s^2=O(s)
\end{equation}
therefore bound all terms containing \(h_s\) and a nonlinear
derivative of \(s\). For the remainder \(\mathcal R\), such a
derivative contains \(l_s^2\) or \(l_v^2\), and the exponential
decay absorbs the factor \(s\) in the last bound.
In terms containing higher derivatives of \(h\), each differentiated
factor of \(s\) is bounded by \(Cl_s\), using \(l_v^2\le l_s\);
\eqref{eq:exponent} then applies. The variable \(x\) is linear
in \(\zeta+\bar\zeta\) and has no higher derivatives.
It follows that, for each \(m\ge2\),
\begin{equation}\label{eq:scaled-potential-bounds}
 \sup_{3B}|D_\eta^m(K_0\circ\Phi_p)|\le C_m.
\end{equation}
Here \(D_\eta^m\) denotes any real coordinate derivative of
order \(m\); the additive constant in \(K_0-K_\infty\) has no
effect.

These estimates hold throughout the chart because
\[
 |x-x(p)|\le C l_x=o(x(p)),\qquad
 |s-s(p)|\le C(l_s+l_v^2)=o(s(p)).
\]
For sufficiently small \(\kappa\), the first-derivative bounds for
the metric extend positivity from the centre to \(3B\).
On \(s\asymp x^{1/4}\), scaled derivatives of the cutoff
contain only products of the bounded factors
\[
 \frac{l_x}{x}=x^{-5/6},\qquad
 \frac{l_s}{s}=x^{-1/12}s^{-3/4},\qquad
 \frac{l_v^2}{s}=x^{-1/6}s^{-3/2}.
\]

\medskip
For \(s(p)\le2S\), we use smooth fibre coordinates on
\(Q(y)=1\). The set \(\{|y|^2\le\cosh(2S)\}\) is compact,
so finitely many coordinate neighbourhoods with fixed domain
margins suffice. We scale both fibre coordinates by \(x^{-1/12}\)
and \(\zeta\) by \(x^{1/6}\), again with the factor \(\kappa\).
For large \(x\), these charts lie where \(K_\infty=K_i\).
Lemma~\ref{lem:inner} gives uniform ellipticity. Smoothness of
\(\Phi\) and \(\Psi\) in \(\tau\) bounds a potential derivative
with \(i\) base and \(j\ge1\) fibre derivatives by
\[
 C_{ij}x^{1/6-5i/6-j/12}\le C_{ij}\qquad(i+j\ge2).
\]
Together with the pure base estimates, this proves
\eqref{eq:scaled-potential-bounds} and uniform metric bounds,
including at \(s=0\).

\medskip
We fix \(\eps_1>0\) small relative to \(\eps_0\). On the central subcone
\[
 |\xi|\le(2-\eps_1)\bigl(\log(2t)-\arsinh1\bigr),
\]
we use the holomorphic coordinates
\begin{equation}\label{eq:central-chart}
 \begin{gathered}
 u=e^\zeta,\qquad z_3=e^\lambda v,\qquad z_1-\ii z_2=e^\lambda,\\
 z_1+\ii z_2=Pe^{-\lambda}-e^\lambda v^2.
 \end{gathered}
\end{equation}
After an \(\SO(3)\) rotation, the centre satisfies \(v=0\) and
\[
 e^{\lambda_p}=\sqrt{t(p)+\sqrt{t(p)^2-|P(u(p))|^2}}>0.
\]
We take steps of size \((\log t(p))^{1/6}\) in both \(\zeta\)
and \(\lambda\), and \((\log t(p))^{-1/3}\) in \(v\).
Homogeneity and Lemma~\ref{lem:F} give, on a slightly larger subcone,
\[
 \frac{c}{(\log t)^{1/3}}I\le D^2F
 \le\frac{C}{(\log t)^{1/3}}I,\qquad
 F_r\asymp(\log t)^{2/3}.
\]
These scales normalize the logarithmic and transverse metric
blocks, respectively. Scaled logarithmic derivatives of order
\(j\ge2\) satisfy
\[
 |D_{\mathrm{scaled}}^jF|
 \le C_j(\log t)^{5/3-5j/6}\le C_j,
 \qquad F_r(\log t)^{-2/3}\asymp1.
\]
For the nonlinear terms, the exact formula is
\[
 t=\frac{e^{\lambda+\bar\lambda}}2
       \left(1+2|v|^2+|Pe^{-2\lambda}-v^2|^2\right).
\]
Thus \(\log(2t)=\lambda+\bar\lambda+2\log(1+|v|^2)\)
up to an exponentially small remainder. Every fixed-order scaled
derivative of \(Pe^{-2\lambda}\) is bounded by
\(\operatorname{poly}(\log t)e^{-c\log t+C(\log t)^{1/6}}\).
As in \eqref{eq:nonlinear}, a nonlinear transverse derivative
contributes \((\log t(p))^{-2/3}\), while
\(F_r=O((\log t(p))^{2/3})\). Their product is bounded.
The bounds for scaled logarithmic derivatives control the remaining
terms. This proves \eqref{eq:scaled-potential-bounds}
and uniform metric bounds on \(3B\).
The coordinates \eqref{eq:central-chart} remain holomorphic at
\(P=0\), and
Lemma~\ref{lem:F} gives smoothness across \(\xi=0\).

\smallskip
To check that these charts cover \(X\), we use the identity
on either base end:
\[
 \log(2t)-\arsinh1
 =x/2+s-\arsinh1+O(e^{-x})+O(e^{-2s}).
\]
For \(s\ge\eps_0x\), this gives
\[
 \frac{|\xi|}{r}\le\frac1{1/2+\eps_0}+o(1)<2-\eps_1
\]
after decreasing \(\eps_1\). These charts overlap those with
\(S\le s\le2\eps_0x\), which overlap the charts with
\(s\le2S\). Points with bounded \(x\) and large \(t\) also lie
in the central subcone. Finitely many ordinary holomorphic charts
cover the remaining compact set. In every chart,
\begin{equation}\label{eq:chart-metric-bounds}
 C^{-1}g_{\mathrm{Euc}}\le\Phi_p^*g_0\le Cg_{\mathrm{Euc}},\qquad
 \|\Phi_p^*g_0\|_{C^j(3B)}\le C_j
 \quad(j\ge0).
\end{equation}
The bounds of order \(j+1\) imply the required \(C^{j,\alpha}(2B)\)
bounds for every \(0<\alpha<1\).

\smallskip
For each \(\eta_0\in2B\), the Euclidean ball \(B(\eta_0,1/2)\) lies
in \(3B\). On these balls, \eqref{eq:chart-metric-bounds} gives
uniform ellipticity and uniform \(C^2\) bounds, hence uniform
bounds for the Christoffel symbols and their first derivatives.
The geodesic equation and its linearization imply that, for
sufficiently small initial velocities, the exponential map at
\(\eta_0\) has differential uniformly close to the identity and is
injective. The domain margin keeps these geodesics inside \(3B\).
Consequently, every point of \(2B\) has a normal neighbourhood
of uniformly positive radius:
\[
 \operatorname{inj}_{(3B,\Phi_p^*g_0)}(\eta_0)\ge\iota_0>0
 \qquad(\eta_0\in2B),
\]
where \(\iota_0\) is independent of \(p\) and \(\eta_0\), and the
exponential maps use the pullback metric on the larger chart.
For the image condition in Hein's and Y.~Li's definitions, we
lift a path starting at \(p\) through \(\Phi_p\), with its lift
starting at \(0\). Until the lift reaches \(\partial(B/2)\),
its length satisfies
\[
 L_{g_0}(\gamma)=L_{\Phi_p^*g_0}(\widetilde\gamma)
 \ge C^{-1/2}L_{g_{\mathrm{Euc}}}(\widetilde\gamma).
\]
Thus any path of length less than \(r_0<1/(2\sqrt C)\) has a
lift contained in \(B/2\), where it extends along the whole path.
Consequently,
\[
 B_{g_0}(p,r_0)\subset\Phi_p(B/2)\subset\Phi_p(B),
\]
with \(r_0\) independent of \(p\).

\medskip
Since \(f=-\log(\omega_0^3/\nu)\) and
\(\omega_0^3/\nu\to1\), \eqref{eq:volumeerror} gives the
same weighted decay for \(f\).
On \(S\le s\le2\eps_0x\), passing to the scaled coordinates
introduces the bounded factors \(l_x/x\), \(l_s/s\), and
\(l_v^2/s\). Base-phase derivatives can introduce powers of
\(l_x\), but the phase-dependent terms are exponentially small
in \(x\). For \(s\le2S\), the weighted estimates are in the
smooth fibre coordinates used above. In the central charts,
logarithmic derivatives contribute \((\log t)^{-5/6}\), and
the transverse step is \((\log t)^{-1/3}\); the phase-dependent
terms decay exponentially in \(\log t\). Since
\(\log t\asymp\log t(p)\) throughout each chart, the chain rule gives
\[
 \|f\circ\Phi_p\|_{C^j(2B)}
 \le\frac{C_j}{(1+\log t(p))^{7/4}}
 \qquad(j\ge0)
\]
for all sufficiently distant centres. Restriction to \(B\) proves
\eqref{eq:fchart}; one further derivative gives its H\"older version.
\end{proof}

\medskip
\Needspace{7\baselineskip}
\begin{lemma}\label{lem:radius-hessian}
The background potential satisfies \(|\nabla^2K_0|_{g_0}\le C\).
\end{lemma}
\begin{proof}
At the centre of a scaled chart, we write
\[
 k=K_0\circ\Phi_p,\qquad g_{i\bar j}=k_{i\bar j}.
\]
The K\"ahler connection gives
\[
 (\nabla^2k)_{i\bar j}=g_{i\bar j},\qquad
 (\nabla^2k)_{ij}
 =k_{ij}-\sum_l(\partial_i g_{j\bar l})
                   \sum_m g^{m\bar l}k_m.
\]
By \eqref{eq:scaled-potential-bounds} and
\eqref{eq:chart-metric-bounds}, the metric coefficients and second
derivatives are bounded. We estimate the two factors in each
connection term separately.

\smallskip
The coordinate order is base, logarithmic fibre, transverse fibre
in \eqref{eq:endchart}; base and two smooth fibre coordinates when
\(s\le2S\); and the two logarithmic coordinates followed by the
transverse coordinate in the central subcone. The first derivatives
satisfy \(|k_i|\le C\Lambda_i\), with
\[
\begin{array}{c|ccc}
 &\Lambda_1&\Lambda_2&\Lambda_3\\ \hline
 S\le s\le2\eps_0x
 &x^{5/6}&x^{1/12}s^{3/4}&1\\
 s\le2S&x^{5/6}&x^{1/12}&x^{1/12}\\
 \text{central subcone}
 &(\log t)^{5/6}&(\log t)^{5/6}&1.
\end{array}
\]
We claim that, at each centre,
\begin{equation}\label{eq:hessian-scales}
 |k_{ij}|\le C,\qquad
 |\partial_i g_{j\bar l}|\le\frac{C}{\Lambda_l},\qquad
 \left|\sum_m g^{m\bar l}k_m\right|\le C\Lambda_l.
\end{equation}
The first bound is \eqref{eq:scaled-potential-bounds}.

\smallskip
On \(S\le s\le2\eps_0x\), the first derivatives have orders
\[
 l_x h_x=O(x^{5/6}),\qquad
 l_s h_s=O(x^{1/12}s^{3/4}),
\]
and the transverse first derivative vanishes at \(v=0\).
A base derivative of a metric coefficient contributes
\(l_x/x=\Lambda_1^{-1}\), and a logarithmic fibre derivative
contributes \(l_s/s=\Lambda_2^{-1}\). For example,
\[
 x^{5/3-3}l_x^3=x^{-5/6},\qquad
 x^{1/6-1}s^{-1/2}l_xl_s^2=x^{-5/6},\qquad
 x^{1/6}s^{-3/2}l_s^3=x^{-1/12}s^{-3/4}.
\]
For the transverse metric coefficient, differentiation of
\(h_s l_v^2\) gives the same factors:
\[
 |l_x h_{xs}l_v^2|\le C\Lambda_1^{-1},\qquad
 |l_s h_{ss}l_v^2|\le C\Lambda_2^{-1}.
\]
The nonlinear terms satisfy the same bounds by
\eqref{eq:nonlinear} and the decay of \(\mathcal R\).
Since \(g_{j\bar l}=k_{j\bar l}\), commuting derivatives gives
\(|\partial_i g_{j\bar l}|\le C/\Lambda_l\) for \(l=1,2\).
For \(l=3\), the uniform third-derivative bound suffices.

At \(v=0\), the mixed transverse entries vanish. The base--fibre
entries of the metric and its inverse satisfy
\[
 |g_{1\bar2}|+|g^{1\bar2}|
 \le C(s/x)^{3/4}=C\frac{\Lambda_2}{\Lambda_1},
 \qquad \Lambda_2\le C\Lambda_1.
\]
Together with uniform ellipticity and \(|k_i|\le C\Lambda_i\),
this gives
\[
 \left|\sum_m g^{m\bar1}k_m\right|
 \le C\left(\Lambda_1+\frac{\Lambda_2^2}{\Lambda_1}\right)
 \le C\Lambda_1,\qquad
 \left|\sum_m g^{m\bar2}k_m\right|\le C\Lambda_2.
\]
The third component is zero, proving \eqref{eq:hessian-scales}
in this region.

\smallskip
When \(s\le2S\), the first derivatives have orders \(x^{5/6}\)
in the base and \(x^{1/12}\) in the fibre. The derivative estimates
above give third derivatives of order \(x^{-5/6}\) whenever a
base derivative occurs, and \(x^{-1/12}\) for three fibre
derivatives. The base--fibre metric block has order
\[
 x^{-3/4}=\frac{x^{1/12}}{x^{5/6}}.
\]
Uniform ellipticity gives the same bound for the inverse mixed
block. Combining it with the first-derivative estimates proves
\eqref{eq:hessian-scales} in this region.

\smallskip
In \eqref{eq:central-chart}, both logarithmic first derivatives are
\(O((\log t)^{5/6})\). Three logarithmic derivatives give
\((\log t)^{5/3-15/6}=(\log t)^{-5/6}\); differentiating the
transverse coefficient gives
\[
 (\log t)^{1/6}F_{rr}(\log t)^{-2/3}
 =O((\log t)^{-5/6}),
\]
and the same bound holds with \(F_{\xi r}\).
These are the bounds \(C/\Lambda_l\) for the two logarithmic
directions. The transverse first derivative
and mixed transverse blocks vanish at \(v=0\), so inversion of
the remaining uniformly elliptic block preserves the common scale
\((\log t)^{5/6}\). The cutoff factors in the gluing region and
the exponentially small phase terms preserve all these estimates.
This establishes \eqref{eq:hessian-scales} in every region.

\smallskip
Substitution into the connection formula gives
\[
 |(\nabla^2k)_{ij}|
 \le C+\sum_{l=1}^3\frac{C}{\Lambda_l}\,C\Lambda_l\le C,
 \qquad |(\nabla^2k)_{i\bar j}|\le C.
\]
Uniform ellipticity in \eqref{eq:chart-metric-bounds} gives
\(|\nabla^2K_0|_{g_0}\le C\) outside a compact set.
Smoothness gives the bound on the compact core.
\end{proof}

\end{document}